\documentclass[english, a4paper]{article}

\usepackage[english]{babel}
\usepackage[T1]{fontenc}
\usepackage{authblk}

\usepackage[a4paper,top=3cm,bottom=3cm,left=3cm,right=3cm,marginparwidth=1.75cm]{geometry}

\usepackage[colorlinks=true, allcolors=blue, hypertexnames=false]{hyperref}
\usepackage{empheq}
\usepackage{graphicx}
\usepackage[colorinlistoftodos]{todonotes}
\usepackage{amsmath}
\usepackage{mathtools}
\usepackage{upgreek}
\usepackage{amssymb}
\usepackage{amsfonts}
\usepackage{amsthm}
\usepackage{hyperref}
\usepackage{enumitem}
\usepackage{mathrsfs}
\usepackage{fontenc}
\usepackage{verbatim}
\usepackage{framed}
\usepackage{algorithm}
\usepackage{bbm}
\usepackage{algpseudocode}
\usepackage{appendix}
\usepackage{color}
\usepackage{multirow}

\numberwithin{equation}{section}
\theoremstyle{plain}
\newtheorem{thm}{Theorem}

\newtheorem{cor}{Corollary}
\newtheorem{assu}{Assumption}

\newtheorem*{rem}{Remark}

\theoremstyle{definition}

\theoremstyle{remark}

\newcommand{\transpose}{^{\operatorname{T}}}

\title{Solving stochastic optimal control problem via stochastic maximum principle with deep learning method}
\author[1]{Shaolin Ji}
\author[2]{Shige Peng}
\author[1]{Ying Peng}
\author[2]{Xichuan Zhang}
\affil[1]{Shandong University-Zhongtai Securities Institute for Financial Studies, Shandong University, 250100, China}
\affil[2]{School of Mathematics, Shandong University, 250100, China}
\begin{document}
\maketitle

\begin{abstract}

In this paper, we aim to solve the high dimensional stochastic optimal control problem from the view of the stochastic maximum principle via deep learning.
By introducing the extended Hamiltonian system which is essentially an FBSDE with a maximum condition,
we reformulate the original control problem as a new one.
Three algorithms are proposed to solve the new control problem.
Numerical results for different examples demonstrate the effectiveness of our proposed algorithms,
especially in high dimensional cases. And an important application of this method is to calculate the sub-linear expectations,
which correspond to a kind of fully nonlinear PDEs.

\end{abstract}

\section{Introduction}\label{sec:intro}

It is well known that Pontryagin's maximum principle~\cite{Bismut1972,Bismut1978,Bensoussan1983Stochastic,pontrygin1987} and Bellman's dynamic programming principle~\cite{Bellman1958Dynamic} are two of the most important tools in solving stochastic optimal control problems. Since these two principles were proposed, the stochastic control theory has been widely developed and extended to a variety of complicated situations in sciences and technologies.

There are many numerical methods for solving stochastic optimal control problems, such as the Markov chain approximation method~\cite{Harold_numerical,Quadrat1994Numerical} which approximate the original controlled process by an appropriate controlled Markov chain on a finite state space, the finite-difference approximations~\cite{dong2007the,jakobsen2003on,krylov2005the} and the probabilistic numerical methods based on dynamic programming~\cite{bertsekas1995neuro-dynamic,PardalosApproximate}. However, few of these methods can deal with high-dimensional problems due to the ``curse of dimensionality''. In other words, the computational complexity grows exponentially when the dimension increases.

In recent years, the deep learning method has been developed rapidly and achieved successes in solving high-dimensional problems of many areas~\cite{BengioDL}, such as computer vision, natural language processing, gaming, etc. This poses a possible way to solve the ``curse of dimensionality''.

Recently, the deep learning method demonstrated remarkable performance in solving the stochastic optimal control problems and the backward stochastic differential equations (BSDEs in short), especially for high dimensional cases~\cite{han2016deep,WeinanDLforBSDE,HanPNAS,deeplearning_FBSDE,Peng_FBSDE_numerical,hure2020deep,raissi2018forward-backward}. The main idea is to treat the control as the parameters in deep neural networks (DNNs in short) and to compute the optimal parameters with stochastic gradient descent methods (SGD). Based on this idea, some researchers extended the neural network architectures to solve the stochastic optimal control problems. For example, \cite{pham2018deep_1,pham2018deep_2} proposed deep learning algorithms from the view of dynamic programming for solving the stochastic control problems. \cite{pereira2019neural} proposed two architectures consisting feed-forward and recurrent neural network to calculate a specific nonlinear stochastic control problem through the Hamilton-Jacobi-Bellman (HJB) equation. The readers can also refer to a recent survey paper which present and compare different deep learning algorithms for solving stochastic control problems and non linear PDEs with the application in finance~\cite{germain2021neural}.

In this paper, different from the above mentioned methods, we solve the stochastic optimal control problem from the view of the stochastic maximum principle (SMP in short) via deep learning. We mainly consider the following stochastic optimal control problem which was introduced in \cite{Peng90,Yong_stochastic_control}:
\begin{align}
\inf_{u(\cdot)\in\mathcal{U}_{ad}[0,T]}\mathbb{E} \Big\{ \int_{0}^{T} f(t, x_t, u_t) \mathrm{d} t + h(x_T) \Big\}, \vspace{1ex} \label{eq:sc_object_sec1}\\
\text{s.t. } x_t = x_0 + \int_{0}^{t} b(t,x_s,u_s)\mathrm{d} s + \int_{0}^{t} \sigma(t,x_s,u_s)\mathrm{d} W_s.
\end{align}

As stated in the SMP, any optimal control along with the optimal state trajectory must be the solutions of the Hamiltonian system plus a maximum condition of a function $H(t,x,u,p,q)$ called the Hamiltonian. In our context, the (extended) Hamiltonian system is characterized by the following FBSDE with a maximum condition~\cite{Bismut1972,Bismut1978,Bensoussan1982}:
\begin{equation}\label{Hamiltonian_sec1}
\left\{
\begin{array}{l}
\mathrm{d} x^*_t = b(t,x^*_t,u^*_t)\mathrm{d} t + \sigma(t,x^*_t,u^*_t)\mathrm{d} W_t, \vspace{1ex}\\
\mathrm{d} p^*_t = -H_x(t,x^*_t,u^*_t,p^*_t,q^*_t)\mathrm{d} t + q^*_t\mathrm{d} W_t, \vspace{1ex}\\
x^*(0) = x_0, \qquad p^*_T=-h_x(x^*_T), \vspace{1ex}\\
H(t,x^*_t,u^*_t,p^*_t,q^*_t) = \max_{u\in U} H (t,x^*_t,u,p^*_t,q^*_t),
\end{array}
\right.
\end{equation}
which has only first-order adjoint equations and correspond to the problems with convex control domain. To the best of our knowledge, this is the first work which solves stochastic optimal control problems through the SMP with deep learning. Our framework is also applicable to more complex problems such as that with non-convex control domain where the Hamiltonian system has second-order adjoint equations, and the cases where the state equation is described by a fully coupled FBSDE. More details are shown in Appendix \ref{appendix:non-convex}.

We first reformulate \eqref{Hamiltonian_sec1} as a new optimal control problem,
\begin{align}
    & \qquad \qquad \begin{matrix}
        \inf_{\tilde{p}_0,\{\tilde{q}_t\}_{0\leq t\leq T}} \mathbb{E}\Big[|-h_x(\tilde{x}_T)-\tilde{p}_T|^2 \Big]
    \end{matrix}\label{eq:new_sc_contr_sec1}\\
&\begin{array}{l}
\mbox{s.t. } \tilde{x}_t = x_0 + \displaystyle \int_{0}^{t} b(t,\tilde{x}_s,\tilde{u}_s)\mathrm{d} s + \int_{0}^{t} \sigma(t,\tilde{x}_s,\tilde{u}_s)\mathrm{d} W_s, \vspace{1ex} \\
\hspace{1.8em} \tilde{p}_t = \tilde{p}_0 - \displaystyle \int_{0}^{t} H_x(s,\tilde{x}_s,\tilde{u}_s,\tilde{p}_s,\tilde{q}_s) \mathrm{d} s + \int_{0}^{t} \tilde{q}_s\mathrm{d} W_s, \vspace{1ex} \\
\hspace{1.8em} \tilde{u}_t = \underset{u\in U}{\arg\max} H (t,\tilde{x}_t,u,\tilde{p}_t,\tilde{q}_t),
\end{array}\notag
\end{align}
where the process $\{\tilde{q}_t\}_{0\leq t\leq T}$ and initial state $\tilde{p}_0$ are regarded as controls. Comparing with \eqref{eq:sc_object_sec1}, the new control problem \eqref{eq:new_sc_contr_sec1} has a simpler quadratic cost functional at time $T$ which provides an easier way to decide whether the state-control pair $(\tilde{x}(\cdot), \tilde{u}(\cdot))$ is an optimal pair, that is whether $\mathbb{E}\Big[|-h_x(\tilde{x}_T)-\tilde{p}_T|^2 \Big]$ equals to 0. And the cost of doing this transformation is that we must deal with an extra term, the maximum condition. In order to solve the new control problem \eqref{eq:new_sc_contr_sec1}, we propose three algorithms suitable for different situations via deep learning. And an important application of our proposed methods is that they can be used to calculate the sub-linear expectations, which correspond to a kind of fully nonlinear PDEs.

In the first Algorithm (Algorithm \ref{alg:1}), a single DNN is constructed to simulate the control $\tilde{q}_t$ and the time $ t $ is regarded as a part of inputs of the neural network. We obtain the approximate estimation of $\tilde{q}_t$ by training such a neural network, and then get the approximate solution $(\tilde{x}_t,\tilde{p}_t,\tilde{q}_t,\tilde{u}_t)_{0\leq t\leq T}$ of \eqref{eq:new_sc_contr_sec1}. For calculating the maximum condition in \eqref{eq:new_sc_contr_sec1}, we employ L-BFGS \cite{LMFGS1989} to approximate the optimal control $\tilde{u}$.

For a general kind of stochastic optimal control problem~\cite{Bismut1972,Bismut1978,Bensoussan1983Stochastic} where all the coefficients are $C^1$ in $u$ and the optimal control $\tilde{u}$ falls inside the control domain, we propose a second algorithm (Algorithm \ref{alg:3}). The aim of this algorithm is to improve the computational efficiency of the approximate solution for the optimal control $\tilde{u}$ in the maximum condition. We first transfer the maximum condition to another kind of constraint $H_u(t,x,u,p,q)=0$, then two neural networks are constructed to simulate the two controls $\{\tilde{q}_t\}_{0\leq t\leq T}$ and $\{\tilde{u}_t\}_{0\leq t\leq T}$, respectively. Moreover, the integral of the constraint $H_u(t,x,u,p,q)$ from 0 to $T$ is added as a penalty term to the original loss function. This algorithm will greatly save the computing time, especially for high dimensional cases  where $\tilde{u}$ can not be solved explicitly.

When the function $\bar{H}$ defined by \eqref{eq:H_bar_definition} is known, we can also solve a class of high-dimensional stochastic optimal control problems even though the optimal control $\tilde{u}$ does not have an explicit solution. For this case, we propose Algorithm 3. Note that when $\tilde{u}$ has an explicit solution as is in the case of Algorithm \ref{alg:1}, we can also get the function $\bar{H}$ explicitly, therefore Algorithm \ref{alg:1} with explicit representation of $\tilde{u}$ is essentially a special case of Algorithm 3.

The numerical results of all the three algorithms demonstrate rather optimistic performance. When $\tilde{u}$ has an explicit solution and thus the function $\bar{H}$ can be solved, Algorithm 3 is an intuitive and better choice. On the other side, even if the optimal control $\tilde{u}$ may not be solved explicitly, our algorithms can still deal with the stochastic optimal control problem. And in this situation, Algorithm \ref{alg:3} or 3 will be better alternatives for high dimensional cases when the conditions mentioned in Section \ref{subsec:numerical_solution 2-NNets} or Section \ref{subsec:legendre} are satisfied. Otherwise Algorithm \ref{alg:1} should be chosen but it is more suitable for low-dimensional cases.

The rest of this paper is structured as follows. In Section \ref{sec:Hamilt_sys}, we briefly introduce the preliminaries about stochastic optimal control problems and reformulate our stochastic optimal control problem as a new control problem. In Section \ref{sec:numerical_alg}, we propose our numerical algorithms for solving the new optimal control problem and present the neural network architecture. In Section \ref{sec:numerical_results}, we show the numerical results and compare the results of our proposed algorithms. More complicated cases with non-convex control domain for solving the second-order adjoint equations are studied in Appendix \ref{appendix:non-convex}.

\section{Preliminaries and problem formulation}\label{sec:Hamilt_sys}

In this section, we introduce the preliminaries of stochastic optimal control and reformulate it to a new control problem.

\subsection{Preliminaries}\label{subsec:Preliminaries}
Let $ T>0 $ and $ (\Omega,\mathcal{F},\mathbb{F},\mathbb{P}) $ be a filtered probability space, where $ W:[0,T] \times \Omega \rightarrow \mathbb{R}^d $ is a $ d $-dimensional standard $ \mathbb{F} $-Brownian motion on $ (\Omega,\mathcal{F},\mathbb{P}) $, $ \mathbb{F}=\{\mathcal{F}_{t}\}_{0\leq t\leq T} $ is the natural filtration generated by the Brownian motion $ W $. Suppose that $ (\Omega,\mathcal{F},\mathbb{P}) $ is complete, $ \mathcal{F}_{0} $ contains all the $ \mathbb{P} $-null sets in $ \mathcal{F} $ and $ \mathbb{F} $ is right continuous. Considering the following \textit{controlled} stochastic differential equation:
\begin{equation}\label{eq:control_system}
  \begin{cases}
    \mathrm{d} x_t = b(t,x_t,u_t)\mathrm{d} t + \sigma(t,x_t,u_t)\mathrm{d} W_t,\\
    x_0 = x_0 \in \mathbb{R}^n,
  \end{cases}
\end{equation}
where $ u_t,t\in[0,T] $, is an admissible control process, i.e. a $ \mathbb{F} $-adapted square-integrable process valued in a given subset $ U $ of $ \mathbb{R}^k $. We define the distance $\|\cdot\|$ in an Euclidean space. $ b $ and $ \sigma $ are the \textit{drift coefficient} and \textit{diffusion coefficient} of \eqref{eq:control_system}, respectively. They are deterministic functions
\begin{align*}
  b: & [0,T] \times \mathbb{R}^n \times U \rightarrow \mathbb{R}^n, \vspace{1ex} \\
  \sigma: & [0,T] \times \mathbb{R}^n \times U \rightarrow \mathbb{R}^{n\times d}.
\end{align*}
The cost functional is
\begin{equation}\label{cost_function}
\begin{matrix}
    J(u(\cdot)) = \mathbb{E} \Big\{ \displaystyle \int_{0}^{T} f(t, x_t, u_t) \mathrm{d} t + h(x_T) \Big\}.
\end{matrix}
\end{equation}

The set of all admissible controls is denoted by $ \mathcal{U}_{ad}[0,T] $
\begin{equation}\label{eq:ad_control_set}
\mathcal{U}_{ad}[0,T] \triangleq \Big\{u:[0,T]\times\Omega\rightarrow U | u \in L_{\mathcal{F}}^2 (0,T;\mathbb{R}^k) \Big\},
\end{equation}
where
\begin{equation*}
L_{\mathcal{F}}^2(0,T;\mathbb{R}^k) \triangleq \left\{x:[0,T]\times\Omega\rightarrow \mathbb{R}^k | x \mbox{ is } \mathbb{F}\mbox{-adapted and } \mathbb{E}[
    \int_{0}^{T}|x_t|^2\mathrm{d} t] <\infty\right\}.
\end{equation*}
Our stochastic optimal control problem can be stated as minimizing \eqref{cost_function} over $ \mathcal{U}_{ad}[0,T] $. The goal is to find $ u^* (\cdot) \in \mathcal{U}_{ad}[0,T]$ (if it exists) such that
\begin{equation}\label{eq:sc_object}
J(u^*(\cdot)) = \inf_{u(\cdot)\in\mathcal{U}_{ad}[0,T]}\mathbb{E} \Big\{
\int_{0}^{T} f(t, x_t, u_t) \mathrm{d} t + h(x_T) \Big\}.
\end{equation}
Any $ u^* (\cdot) \in \mathcal{U}_{ad}[0,T]$ satisfying \eqref{eq:sc_object} is called an \textit{optimal control}. The corresponding state process $ x^*(\cdot) $ and the state-control pair $ (x^*(\cdot), u^*(\cdot)) $ are called an \textit{optimal state process} and an \textit{optimal pair} respectively.

Firstly let us make the following assumptions.
\begin{assu}\label{assu:1}
\begin{enumerate}
  \item[(i)] The maps $ b,\sigma,f $ and $ h $ are measurable, and there exist a constant $ L>0 $ and a modulus of continuity $ \bar{\omega}:\left[0,\infty\right)\rightarrow\left[0,\infty\right) $ such that for $ \varphi(t,x,u)=b(t,x,u),\sigma(t,x,u)$, $f(t,x,u),h(x) $, we have
  \begin{equation}
  \left\{
  \begin{array}{l}
  |\varphi(t,x,u)-\varphi(t,\hat{x},\hat{u})|\leq L|x-\hat{x}| + \bar{\omega}\|u-\hat{u}\|, \\
  \hspace{9em}\forall t\in[0,T], \ \ x,\hat{x}\in\mathbb{R}^{n}, \ \ u,\hat{u}\in U, \\
  |\varphi(t,0,u)| \leq L, \ \ \forall t\in[0,T], u\in U;
  \end{array}
  \right.
  \end{equation}
  \item[(ii)] The maps $ b,\sigma,f $ and $ h $ are $ C^2 $ in $ x $. Moreover, there exist a constant $ L>0 $ and a modulus of continuity $ \bar{\omega}:\left[0,\infty\right)\rightarrow\left[0,\infty\right) $ such that for $ \varphi=b,\sigma,f,h $, we have
  \begin{equation}
  \left\{
  \begin{array}{l}
  |\varphi_{x}(t,x,u)-\varphi_{x}(t,\hat{x},\hat{u})|\leq L|x-\hat{x}| + \bar{\omega}\|u-\hat{u}\|, \\
  |\varphi_{xx}(t,x,u)-\varphi_{xx}(t,\hat{x},\hat{u})|\leq\bar{\omega}(|x-\hat{x}| + \|u-\hat{u}\|), \\
  \hspace{9em}\forall t\in[0,T], \ \ x,\hat{x}\in\mathbb{R}^{n}, \ \ u,\hat{u}\in U.
  \end{array}
  \right.
  \end{equation}
\end{enumerate}
\end{assu}

\begin{assu}\label{assu:2}
The control domain $ U $ is a convex body in $ \mathbb{R}^{k} $. The maps $ b,\sigma $ and $ f $ are locally Lipschitz in $ u $, and their derivatives in $ x $ are continuous in $ (x,u) $.
\end{assu}

In the following, before introducing a set of sufficient conditions for the \textit{Stochastic Maximum Principle} (SMP in short), we firstly introduce the adjoint equations involved in a SMP and the associated stochastic Hamiltonian system.

Let $ (x^*(\cdot),u^*(\cdot)) $ be a given optimal pair. We introduce the adjoint BSDE as follows:
\begin{equation}\label{eq:1-order-adjoint}
\left\{
\begin{array}{l}
\mathrm{d} p_t^* = -\Big\{b_x(t,x^*_t,u^*_t)\transpose p_t^* + \sum_{j=1}^{d}\sigma_x^j(t,x^*_t,u^*_t)\transpose q_{jt}^* - f_x(t,x^*_t,u^*_t) \Big\} \mathrm{d} t + q_t^*\mathrm{d} W_t, \vspace{1ex} \\
p_T^*=-h_x(x^*_T), \qquad t\in[0,T].
\end{array}
\right.
\end{equation}
where $ p^*(\cdot) $ and $ q^*(\cdot) $ are two $ \mathbb{F} $-adapted processes which should be solved. Any pair of processes $ (p^*(\cdot),q^*(\cdot)) $ $\in L_{\mathcal{F}}^2(0,T;\mathbb{R}^n)\times(L_{\mathcal{F}}^2(0,T;\mathbb{R}^n))^d $ satisfying \eqref{eq:1-order-adjoint} is called an \textit{adapted solution} of \eqref{eq:1-order-adjoint}. Under Assumption \ref{assu:1}, for any $ (x^*(\cdot),u^*(\cdot))\in L_{\mathcal{F}}^2(0,T;\mathbb{R}^n) \times \mathcal{U}[0,T] $, \eqref{eq:1-order-adjoint} admits a \textit{unique} adapted solution $ (p^*(\cdot),q^*(\cdot)). $

We refer to \eqref{eq:1-order-adjoint} as the \textit{first-order} \textit{adjoint equations} and to $ p^*(\cdot) $ as the \textit{first-order} \textit{adjoint process}. If $ (x^*(\cdot),u^*(\cdot)) $ is an optimal (resp. admissible) pair, and $ (p^*(\cdot),q^*(\cdot)) $ is an adapted solution of \eqref{eq:1-order-adjoint}, then $ (x^*(\cdot), u^*(\cdot), p^*(\cdot), q^*(\cdot)) $ is called an \textit{optimal 4-tuple} (resp. \textit{admissible 4-tuple}). According to Theorem 5.2 of Chapter 3 and the comments after it in \cite{Yong_stochastic_control}, we have the following sufficient conditions for the SMP:
\begin{thm}\label{thm:suff_cond}
  Let Assumptions \ref{assu:1} and \ref{assu:2} hold. Let $ (x^*(\cdot), u^*(\cdot), $ $ p^*(\cdot), q^*(\cdot)) $ be an admissible 4-tuple. Suppose that $ h(\cdot) $ is convex, $ H(t,\cdot,\cdot,p_t^*,q_t^*) $ defined by
  \begin{equation}\label{eq:Ham_con}
  \begin{array}{l}
  H(t,x,u,p,q) = \left\langle  p,b(t,x,u) \right\rangle + \mbox{tr}[q\transpose\sigma(t,x,u)] - f(t,x,u), \vspace{1ex} \\
  \hspace{7em}(t,x,u,p,q)\in[0,T]\times\mathbb{R}^n\times U\times \mathbb{R}^n \times \mathbb{R}^{n\times d},
  \end{array}
  \end{equation}
  is concave for all $ t\in[0,T] $ almost surely and
  \begin{equation}\label{eq:max_cond}
  H(t,x^*_t,u^*_t,p^*_t,q^*_t) = \max_{u\in U}H(t,x^*_t,u,p^*_t,q^*_t), \qquad \mbox{a.e. }t\in[0,T], \qquad \mathbb{P}\mbox{-a.s.}
  \end{equation}
  holds. Then $ (x^*(\cdot), u^*(\cdot)) $ is an optimal pair of \eqref{eq:sc_object}.
\end{thm}

Note that the partial differentials of the Hamiltonian $H$ satisfy $b(t,x,u)=H_p(t,x,u,p,q)$ and $\sigma(t,x,u)=H_q(t,x,u,p,q)$, then the combination of \eqref{eq:control_system}, \eqref{eq:1-order-adjoint} and \eqref{eq:max_cond} can be written as follows:
\begin{equation}\label{eq:FBSDE_sys}
\left\{
\begin{array}{l}
\mathrm{d} x_t^* = b(t,x^*_t,u^*_t)\mathrm{d} t + \sigma(t,x^*_t,u^*_t)\mathrm{d} W_t, \vspace{1ex} \\
\mathrm{d} p_t^* = -H_x(t,x^*_t,u^*_t,p_t^*,q_t^*)\mathrm{d} t + q_t^*\mathrm{d} W_t, \qquad t\in[0,T], \vspace{1ex} \\
x^*_0 = x_0, \qquad p_T^*=-h_x(x^*_T), \vspace{1ex} \\
H(t,x^*_t,u^*_t,p^*_t,q^*_t) = \max_{u\in U}H(t,x^*_t,u,p^*_t,q^*_t), \vspace{1ex}
\end{array}
\right.
\end{equation}
which is called a \textit{(extended) stochastic Hamiltonian system}, with its solution being a 4-tuple $ (x^*(\cdot), u^*(\cdot), p^*(\cdot), q^*(\cdot)) $. And there exists a function $\bar{H}: [0,T] \times \mathbb{R}^{n} \times \mathbb{R}^{n} \times \mathbb{R}^{n\times d} \rightarrow \mathbb{R}$ such that
\begin{equation}\label{eq:H_bar_definition}
  \bar{H}(t,x,p,q) = \max_{u\in U} H(t,x,u,p,q),
\end{equation}
which is a function independent of control $u$.

In this paper, we primarily study the problems with convex control domain which correspond to the first-order adjoint equations, and more complicated cases are discussed in Appendix \ref{appendix:non-convex}. In order to make sure that the numerical algorithms can completely solve the optimal control problem, we mainly focus on the cases when equation \eqref{eq:control_system} has unique optimal control $u^*$ and its corresponding Hamiltonian system (see \eqref{eq:FBSDE_pure}) has unique adapted solution $(x^*(\cdot),p^*(\cdot),q^*(\cdot))$. These cases do exist when more strictly convex assumptions in Theorem \ref{thm:suff_cond} and the monotonic conditions~\cite{Peng1999Fully} of FBSDE \eqref{eq:FBSDE_pure} hold.

\subsection{Problem formulation}\label{subsec:Problem description}

In this subsection, we reformulate the control problem \eqref{eq:sc_object} as a new problem based on the SMP and its corresponding stochastic Hamiltonian system. Considering the (extended) stochastic Hamiltonian system \eqref{eq:FBSDE_sys}, which is essentially a coupled FBSDE with a maximum condition. Suppose that there exists a solution $ (x^*(\cdot),u^*(\cdot),p^*(\cdot),q^*(\cdot)) $ for FBSDE \eqref{eq:FBSDE_sys}.

As is known, the FBSDE can be regarded as a stochastic optimal control problem~\cite{ma1999forward}. Based on this idea, we have the following state equation with a maximum condition which is equivalent to \eqref{eq:FBSDE_sys},
\begin{equation}\label{eq:for_stoch_diff_eq}
\left\{
\begin{array}{l}
\mathrm{d} \tilde{x}_t = b(t,\tilde{x}_t,\tilde{u}_t)\mathrm{d} t + \sigma(t,\tilde{x}_t,\tilde{u}_t)\mathrm{d} W_t, \vspace{1ex} \\
\mathrm{d} \tilde{p}_t = -H_x(t,\tilde{x}_t,\tilde{u}_t,\tilde{p}_t,\tilde{q}_t)\mathrm{d} t + \tilde{q}_t\mathrm{d} W_t, \vspace{1ex} \\
\tilde{x}_0 = x_0, \qquad \tilde{p}_0=\tilde{p}_0, \vspace{1ex} \\
H(t,\tilde{x}_t,\tilde{u}_t,\tilde{p}_t,\tilde{q}_t) = \max_{u\in U} H (t,\tilde{x}_t,u,\tilde{p}_t,\tilde{p}_t). \vspace{1ex}
\end{array}
\right.
\end{equation}
where $(\tilde{p}_0,\tilde{q})$ is the pair of control valued in $\mathbb{R}^n\times\mathbb{R}^{n\times d}$. As $\tilde{u}_t$ can be represented as
\begin{equation}\label{eq:U_condition}
    \tilde{u}_t = \arg\max_{u\in U}H(t,\tilde{x}_t,u,\tilde{p}_t,\tilde{q}_t),
\end{equation}
then we get a new variational problem which is a reformulation of the control problem \eqref{eq:sc_object}:
\begin{equation}\label{eq:new_sc_contr}
\inf_{\tilde{p}_0,\{\tilde{q}_t\}_{0\leq t\leq T}}\mathbb{E}\Big[|-h_x(\tilde{x}_T)-\tilde{p}_T|^2 \Big]
\end{equation}
\begin{equation*}
\begin{array}{l}
\mbox{s.t. } \tilde{x}_t = x_0 + \displaystyle \int_{0}^{t} b(t,\tilde{x}_s,\tilde{u}_s)\mathrm{d} s + \int_{0}^{t} \sigma(t,\tilde{x}_s,\tilde{u}_s)\mathrm{d} W_s, \vspace{1ex} \\
\hspace{1.8em} \tilde{p}_t = \tilde{p}_0 - \displaystyle \int_{0}^{t} H_x(s,\tilde{x}_s,\tilde{u}_s,\tilde{p}_s,\tilde{q}_s) \mathrm{d} s + \int_{0}^{t} \tilde{q}_s\mathrm{d} W_s, \vspace{1ex} \\
\hspace{1.8em} \tilde{u}_t = \underset{u\in U}{\arg\max} H (t,\tilde{x}_t,u,\tilde{p}_t,\tilde{q}_t).
\end{array}
\end{equation*}

According to the following Theorem \ref{thm:main} and Corollary \ref{cor:equivalent_FBSDE}, we can prove that the optimal control $\tilde{u}$ of \eqref{eq:new_sc_contr} can be obtained when $\mathbb{E}\Big[|-h_x(\tilde{x}_T)-\tilde{p}_T|^2 \Big] = 0$. Supposing that $\bar{H}$ in \eqref{eq:H_bar_definition} is differentiable in $x,p,q$, we have
\begin{equation}\label{eq:notations}
\left\{
\begin{array}{l}
\bar{H}_p(t,x,p,q)=H_p(t,x,u^*,p,q)=b(t,x,u^*), \vspace{1ex} \\
\bar{H}_q(t,x,p,q)=H_q(t,x,u^*,p,q)=\sigma(t,x,u^*), \vspace{1ex} \\
\bar{H}_x(t,x,p,q)=H_x(t,x,u^*,p,q), \vspace{1ex} \\
u^* = \underset{u\in U}{\arg\max} H (t,x,u,p,q),
\end{array}
\right.
\end{equation}
for any $(t,x,u,p,q)\in[0,T]\times\mathbb{R}^n\times U\times \mathbb{R}^n \times \mathbb{R}^{n\times d}$. Then \eqref{eq:FBSDE_sys} can be rewritten as
\begin{equation}\label{eq:FBSDE_pure}
\left\{
\begin{array}{l}
\mathrm{d} x_t^* = \bar{H}_p(t,x_t^*,p_t^*,q_t^*)\mathrm{d} t + \bar{H}_{q}(t,x_t^*,p_t^*,q_t^*)\mathrm{d} W_t, \vspace{1ex} \\
\mathrm{d} p_t^* = -\bar{H}_x(t,x_t^*,p_t^*,q_t^*)\mathrm{d} t + q_t^*\mathrm{d} W_t, \vspace{1ex} \\
x_0^* = x_0, \qquad p_T^* = -h_x(x^*_T),
\end{array}
\right.
\end{equation}
which is a FBSDE without constraint. And the problem \eqref{eq:new_sc_contr} is equivalent to
\begin{equation}\label{eq:new_sc_contr_pure}
\inf_{\tilde{p}_0,\{\tilde{q}_t\}_{0\leq t\leq T}}\mathbb{E}\Big[|-h_x(\tilde{x}_T)-\tilde{p}_T|^2 \Big]
\end{equation}
\begin{equation*}
\begin{array}{l}
\mbox{s.t. } \tilde{x}_t = x_0 + \displaystyle \int_{0}^{t} \bar{H}_p(s,\tilde{x}_s,\tilde{p}_s,\tilde{q}_s)\mathrm{d} s + \int_{0}^{t} \bar{H}_q(s,\tilde{x}_s,\tilde{p}_s,\tilde{q}_s)\mathrm{d} W_s, \vspace{1ex} \\
\hspace{1.8em} \tilde{p}_t = \tilde{p}_0 - \displaystyle \int_{0}^{t} \bar{H}_x(s,\tilde{x}_s,\tilde{p}_s,\tilde{q}_s) \mathrm{d} s + \int_{0}^{t} \tilde{q}_s\mathrm{d} W_s.
\end{array}
\end{equation*}

Define the cost functional of the stochastic optimal control problem \eqref{eq:new_sc_contr_pure} as
\begin{equation*}{}
    V(0,x_0) = \inf_{\tilde{p}_0,\{\tilde{q}_t\}_{0\leq t\leq T}}\mathbb{E}\Big[|-h_x(\tilde{x}_T)-\tilde{p}_T|^2 \Big],
\end{equation*}
and we have the following theorem which can be referred to Proposition 1.1 in Chapter 3 of \cite{ma1999forward}.
\begin{thm}\label{thm:main}
    Assume $h_x(x)$ is continuous in $x$, the map $\bar{H}$ defined by \eqref{eq:H_bar_definition} is differentiable in $x,p,q$, and its derivatives are continuous in $x,p,q$, such that for $\varphi=\bar{H}_x, \bar{H}_p, \bar{H}_q$ and some constant $C>0$,
    \begin{align*}
        |\varphi(t,x,p,q)-\varphi(t,x',p',q')| &\leq C(|x-x'|+|p-p'|+|q-q'|), \\
        |\varphi(t,0,0,0)|,|\bar{H}_q(t,x,p,0)|&\leq C,\\
        & \forall x,x',p,p' \in \mathbb{R}^n, q,q'\in \mathbb{R}^{n\times d},
    \end{align*}
    then the FBSDE \eqref{eq:FBSDE_pure} is solvable over $[0,T]$ if and only if $V(0,x_0)=0$.
\end{thm}

\begin{cor}\label{cor:equivalent_FBSDE}
    Suppose the assumptions in Theorem \ref{thm:suff_cond}, \ref{thm:main} hold, if there exists a solution $ (\tilde{x}_t, \tilde{u}_t, \tilde{p}_t, \tilde{q}_t)_{0\leq t\leq T} $ of \eqref{eq:for_stoch_diff_eq} satisfying
    \begin{equation}\label{eq:obj_fun}
    \mathbb{E}\Big[|-h_x(\tilde{x}_T)-\tilde{p}_T|^2 \Big] = 0,
    \end{equation}
    then $(\tilde{x}_t, \tilde{u}_t, \tilde{p}_t, \tilde{q}_t)$ is a solution of \eqref{eq:FBSDE_sys}, and the cost functional \eqref{eq:sc_object} can be obtained by
    \begin{equation}
    \begin{matrix}
        J(u^*(\cdot))=J(\tilde{u}(\cdot)) = \mathbb{E}\Big\{ \displaystyle \int_{0}^{T} f(t, \tilde{x}_t, \tilde{u}_t) \mathrm{d} t + h(\tilde{x}_T) \Big\}.
    \end{matrix}
    \end{equation}
\end{cor}

\begin{rem}\label{rem:1}
In this corollary, the assumptions in Theorem \ref{thm:suff_cond} ensure that the optimal pair $(x^*(\cdot),u^*(\cdot))$ can be got by solving the Hamiltonian system \eqref{eq:FBSDE_sys}, the assumptions in Theorem \ref{thm:main} ensure that the optimal control of $\tilde{u}$ in \eqref{eq:new_sc_contr} exists so that we can solve the Hamiltonian system \eqref{eq:FBSDE_sys} by solving the optimal control problem \eqref{eq:new_sc_contr}.
\end{rem}

Problem \eqref{eq:new_sc_contr} is a reformulation of the problem \eqref{eq:FBSDE_sys}. This new control problem has a quadratic cost functional at the terminal time, and it provides an alternative criterion to decide whether the state-control pair $(\tilde{x}(\cdot), \tilde{u}(\cdot))$ is an optimal pair, which is whether $\mathbb{E}\Big[|-h_x(\tilde{x}_T)-\tilde{p}_T|^2 \Big]$ equals to 0. But the cost is that we must deal with an additional maximum condition $\tilde{u}_t = \arg\max_{u\in U} H (t,\tilde{x}_t,u,\tilde{p}_t,\tilde{q}_t)$. And in the rest of this paper, we mainly focus on solving problem \eqref{eq:new_sc_contr} through the deep learning method.

\section{Numerical algorithms}\label{sec:numerical_alg}

In Section \ref{sec:Hamilt_sys}, we briefly introduce the sufficient conditions of the SMP and reformulate the Hamiltonian system to a new variational problem. In this section, we propose three algorithms for solving the new variational problem through deep learning and show the neural network structure of Algorithm \ref{alg:1}. The network structures of the other two algorithms are similar with that of Algorithm \ref{alg:1}.

\subsection{Algorithm 1: Numerical algorithm with 1-NNet}\label{subsec:numerical_solution}

Let $\pi$ be a partition of the time interval, $ 0 = t_{0}<t_{1}<t_{2}<\cdots<t_{N-1}<t_{N} =T $ of $ [0,T] $. Define $ \Delta t_{i}=t_{i+1}-t_{i} $ and $ \Delta W_{t_i}=W_{t_{i+1}}-W_{t_{i}} $, where $ W_{t_{i}} \sim \mathcal{N}(0,t_{i}) $, for $ i = 0, 1, 2,\cdots, N-1 $. We also denote
$$ \delta = \sup_{0\leq i\leq N-1}\Delta t_i. $$
Then the Euler scheme of the forward SDE \eqref{eq:for_stoch_diff_eq} and the maximum condition \eqref{eq:U_condition} can be written as
\begin{equation}\label{eq:for_stoch_diff_eq1}
\left\{
\begin{array}{l}
\tilde{x}^{\pi}_{t_{i+1}} = \tilde{x}^{\pi}_{t_i} + b(t_i,\tilde{x}^{\pi}_{t_i},\tilde{u}^{\pi}_{t_i})\Delta t_i + \sigma(t_i,\tilde{x}^{\pi}_{t_i},\tilde{u}^{\pi}_{t_i})\Delta W_{t_i}, \vspace{1ex} \\
\tilde{p}^{\pi}_{t_{i+1}} = \tilde{p}^{\pi}_{t_i}-H_x(t_i,\tilde{x}^{\pi}_{t_i},\tilde{u}^{\pi}_{t_i},\tilde{p}^{\pi}_{t_i},\tilde{q}^{\pi}_{t_i})\Delta t_i + \tilde{q}^{\pi}_{t_i}\Delta W_{t_i}, \vspace{1ex} \\
\tilde{x}^{\pi}_0 = x_0, \qquad \tilde{p}^{\pi}_0=\tilde{p}_0, \vspace{1ex} \\
\tilde{u}^{\pi}_{t_i} = \underset{u\in U}{\arg\max} H (t_i,\tilde{x}^{\pi}_{t_i},u, \tilde{p}^{\pi}_{t_i}, \tilde{q}^{\pi}_{t_i}).
\end{array}
\right.
\end{equation}

We regard $ \{\tilde{q}^{\pi}_{t_i}\}_{0\leq i<N} $ as a control and assume it satisfying
\begin{equation}\label{eq:pro_fun_q}
\tilde{q}^{\pi}_{t_i} = \phi^1(t_i, \tilde{x}^{\pi}_{t_i},\tilde{u}^{\pi}_{t_i},\tilde{p}^{\pi}_{t_i};\theta^1_{t_i}),
\end{equation}
where $\tilde{q}^{\pi}_{t_i}$ is a feedback control of the states $ \tilde{x}^{\pi}_{t_i}$, $\tilde{p}^{\pi}_{t_i} $ and the control $ \tilde{u}^{\pi}_{t_i} $.
 Here we give some explications to recognize the difference between the controls $ \tilde{u}^{\pi}_{t_i} $ and $\tilde{q}^{\pi}_{t_i}$. $ \tilde{u}^{\pi}_{t_i} $ is the control of the original control problem \eqref{eq:control_system} where it should be solved by the maximum condition \eqref{eq:U_condition} in \eqref{eq:for_stoch_diff_eq}. And $ \tilde{q}^{\pi}_{t_i} $ is the control of the new control problem \eqref{eq:new_sc_contr}.

Note that the processes $\tilde{q}_{t_{i}}^{\pi}$ and $\tilde{u}_{t_{i}}^{\pi}$ are interdependent according to \eqref{eq:U_condition} and \eqref{eq:pro_fun_q}. Plugging \eqref{eq:U_condition}  in \eqref{eq:pro_fun_q} and according to the implicit function theorem, we can get
\begin{align}\label{eq:pro_fun_q_new}
\tilde{q}^{\pi}_{t_i} &= \phi(t_i, \tilde{x}^{\pi}_{t_i},\tilde{p}^{\pi}_{t_i};\theta^1_{t_i}) \vspace{1ex} \\
&=\phi^1(t_i, \tilde{x}^{\pi}_{t_i},\arg\max_{u\in U} H (t_i,\tilde{x}^{\pi}_{t_i},u, \tilde{p}^{\pi}_{t_i}, \tilde{q}^{\pi}_{t_i}),\tilde{p}^{\pi}_{t_i};\theta^1_{t_i}), \notag
\end{align}
where $ \phi $ is a new unknown function.

\begin{figure}[H]
  \centering
  \includegraphics[scale=0.8]{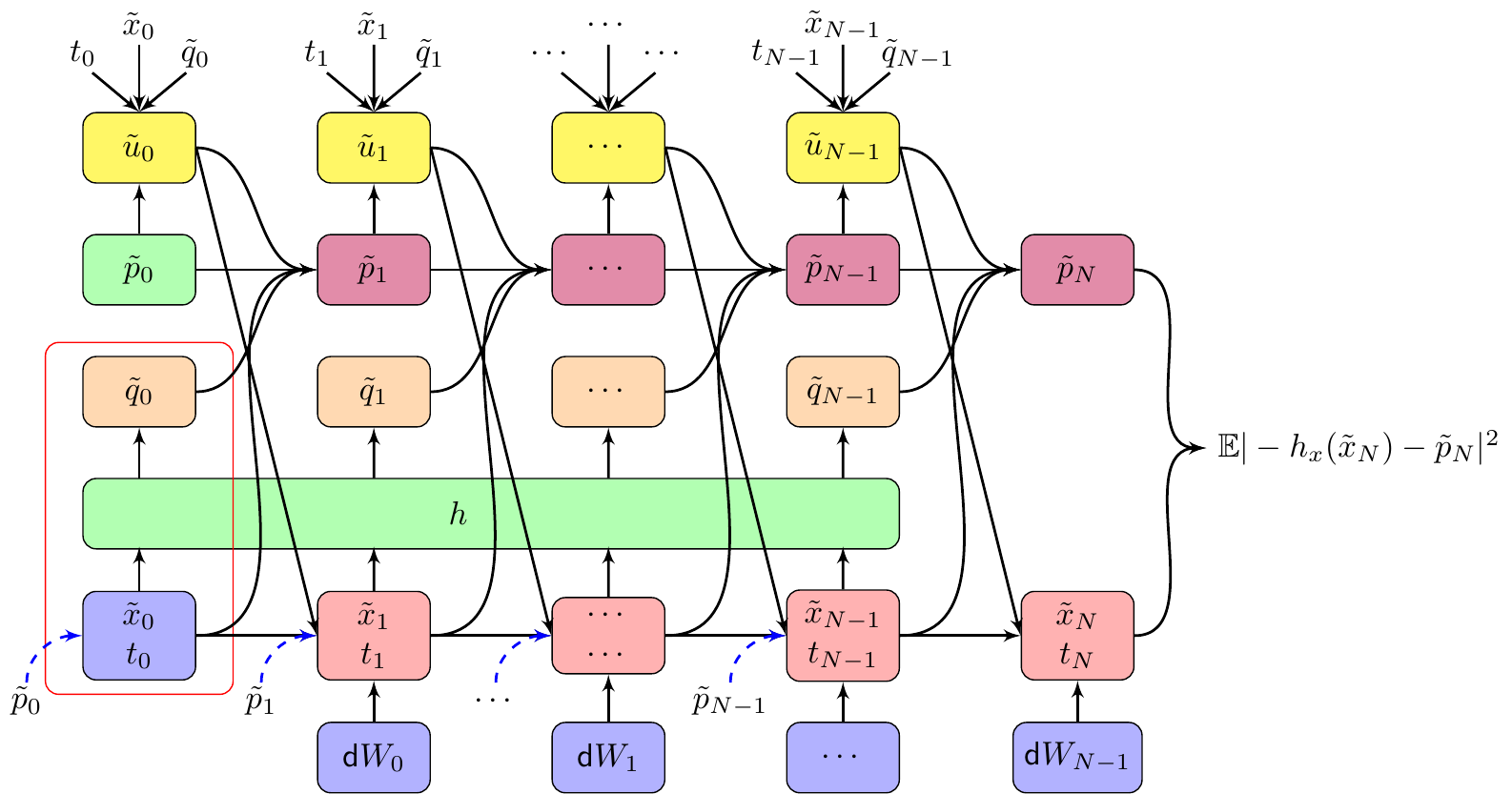}
  \caption{The whole network architecture for Algorithm \ref{alg:1}. $h$ represents the hidden layers which are the common for different time-points and $\tilde{q}_i$ represents the outputs. The value of $\tilde{q}_0$, the weights and biases of the hidden layers are trainable parameters. $u_i$ is a function of $(t_i,\tilde{x}_i,\tilde{p}_i,\tilde{q}_i)$. The box with red line represents the DNN of a single time-point $t_0$. The solid lines represent the data flow generated in the current iteration, the dashed blue lines represent the neural network inputs at each time-point. }
  \label{fig:whole_net}
\end{figure}

We develop a neural network (1-NNet) for simulating the feedback control $ \tilde{q}^{\pi}_{\cdot} $. Different from our previous work in \cite{Peng_FBSDE_numerical}, a single network is constructed for all the time-points and the time $ t_{i} $ is regarded as an input of the neural network. The network consists of five layers including one $ (1+n+n) $-dim input layer, three $ (10+n+n) $-dim hidden layers and a $ (n\times d) $-dim output layer. All parameters of the network are represented as $ \theta $. The loss function is defined as
\begin{equation}\label{eq:cost_fun_1}
\mbox{loss} = \dfrac{1}{M} \sum_{j=1}^{M}\Big[ |-h_x(\tilde{x}^{\pi}_T)-\tilde{p}^{\pi}_T|^2 \Big],
\end{equation}
where $ M $ is the number of samples. Figure~\ref{fig:whole_net} gives the whole network structure for all the time-points. Moreover, we use a box with red line in Figure~\ref{fig:whole_net} to show the DNN of a single time-point $t_0$.

For convenience, the time interval $ [0, T] $ is partitioned evenly, i.e. $ \Delta t_i = t_{i+1}-t_i=T/N $ for all $ i=0,1,\cdots,N,N\geq 1 $. We define $ \Delta W_{t_i}=W_{t_{i+1}}-W_{t_{i}} $ and denote the iteration step by $ l $ which is marked by superscript in the algorithm. The pseudo-code for solving the stochastic optimal control problem is given in Algorithm \ref{alg:1}.

\begin{algorithm}[H]
  \renewcommand{\thealgorithm}{1}
  \caption{Numerical algorithm with 1-NNet}
  \label{alg:1}
  \begin{algorithmic}[1]
    \Require The Brownian motion $ \Delta W_{t_i} $, initial parameters $ (\theta^0,\tilde{p}_0^{0,\pi}) $, learning rate $ \eta $;
    \Ensure The 4-tuple precesses $ (\tilde{x}^{l,\pi}_{t_i},\tilde{u}^{l,\pi}_{t_i},\tilde{p}^{l,\pi}_{t_i},\tilde{q}^{l,\pi}_{t_i}) $.
    \For { $ l = 0 $ to $ maxstep $}
    \State $ \tilde{x}_{0}^{l,\pi} = x_{0} $, $ \tilde{p}_{0}^{l,\pi} = \tilde{p}_{0}^{l,\pi}; $
    \For { $ i = 0 $ to $ N-1 $}
    \State $\tilde{q}^{l,\pi}_{t_i} = \phi(t_i, \tilde{x}^{l,\pi}_{t_i},\tilde{p}^{l,\pi}_{t_i};\theta^{l});$
    \State $\tilde{u}^{l,\pi}_{t_i} = \arg\max_{u\in U}H(t_i,\tilde{x}^{l,\pi}_{t_i},u,\tilde{p}^{l,\pi}_{t_i},\tilde{q}^{l,\pi}_{t_i});$
    \State $\tilde{x}^{l,\pi}_{t_{i+1}} = \tilde{x}^{l,\pi}_{t_i} + b(t_i,\tilde{x}^{l,\pi}_{t_i},\tilde{u}^{l,\pi}_{t_i})\Delta t_i + \sigma(t_i,\tilde{x}^{l,\pi}_{t_i},\tilde{u}^{l,\pi}_{t_i})\Delta W_{t_i};$
    \State $\tilde{p}^{l,\pi}_{t_{i+1}} = \tilde{p}^{l,\pi}_{t_i}-H_x(t_i,\tilde{x}^{l,\pi}_{t_i},\tilde{u}^{l,\pi}_{t_i},\tilde{p}^{l,\pi}_{t_i},\tilde{q}^{l,\pi}_{t_i})\Delta t_i + \tilde{q}^{l,\pi}_{t_i}\Delta W_{t_i};$
    \EndFor
    \State $ J(\tilde{u}^{l,\pi}(\cdot))=\dfrac{1}{M} \sum_{j=1}^{M}\Big[\dfrac{T}{N} \sum_{i=0}^{N-1}f(t_i,\tilde{x}_{t_i}^{l,\pi},\tilde{u}_{t_i}^{l,\pi}) + h(\tilde{x}_T^{l,\pi}) \Big]; $
    \State $ \mbox{loss} = \dfrac{1}{M} \sum_{j=1}^{M}\Big[ |-h_x(\tilde{x}^{l,\pi}_T)-\tilde{p}^{l,\pi}_T|^2 \Big]; $
    \State $ (\theta^{l+1},\tilde{p}_0^{l+1,\pi})=(\theta^{l},\tilde{p}_0^{l,\pi})-\eta\nabla\mbox{loss}. $
    \EndFor
  \end{algorithmic}
\end{algorithm}

\begin{rem}\label{rem:2}
    There is one difficulty in this algorithm. As shown in line 5 of Algorithm \ref{alg:1}, the solution of the maximum condition is an extremum problem of multivariate functions and  in most cases it has no analytical solution, which means that we can not get the explicit value of $ \tilde{u}^{l,\pi}_{t_i} $ to calculate the forward process. When the explicit solution is not available, there are some ways to get the approximated solution, such as the BFGS and its extended methods \cite{LMFGS1989}, the gradient descent methods, the Sequential Least Squares Programming (SLSQP) and so on. On the other hand, when the explicit solution of $\tilde{u}$ is available, the Hamiltonian system \eqref{eq:for_stoch_diff_eq} is equivalent to a FBSDE, and this situation will be discussed in subsection \ref{subsec:legendre}.
\end{rem}

\subsection{Algorithm 2: Numerical algorithm with 2-NNets}\label{subsec:numerical_solution 2-NNets}

As mentioned in Algorithms \ref{alg:1}, when the explicit solution of the maximum condition is not available, some approximation methods should be used. However, it is very time consuming to calculate the approximate maximum condition for high-dimensional cases. In order to solve this problem, we develop a numerical algorithm with two neural networks (2-NNets) in this subsection. Here we consider a kind of stochastic optimal control problem, where the convex control domain $ U=\mathbb{R}^k $, all the coefficients are $ C^1 $ in $ u $ and the optimal control $\tilde{u}$ falls inside the boundary of the control domain, then \eqref{eq:max_cond} implies
\begin{equation}
\begin{array}{l}
H_u(t,x^*_t,u^*_t,p_t^*,q_t^*) = 0, \qquad \forall u\in U, \qquad\mbox{a.e. }t\in[0,T], \qquad \mathbb{P}\mbox{-a.s.}
\end{array}
\end{equation}
Thus the corresponding stochastic Hamiltonian system \eqref{eq:FBSDE_sys} can be represented as
\begin{equation}\label{eq:Hamil_Case_Hu}
\left\{
\begin{array}{l}
\mathrm{d} x^*_t = b(t,x^*_t,u^*_t)\mathrm{d} t + \sigma(t,x^*_t,u^*_t)\mathrm{d} W_t, \vspace{1ex} \\
\mathrm{d} p_t^* = -H_x(t,x^*_t,u^*_t,p_t^*,q_t^*)\mathrm{d} t + q_t^*\mathrm{d} W_t, \qquad t\in[0,T], \vspace{1ex} \\
x^*_0 = x_0, \qquad p_T^*=-h_x(x^*_T), \vspace{1ex} \\
H_u(t,x^*_t,u^*_t,p_t^*,q_t^*) = 0, \qquad \forall u\in U.
\end{array}
\right.
\end{equation}

The Hamiltonian system \eqref{eq:Hamil_Case_Hu} with constraint $H_u(t,x^*_t,u^*_t,p_t^*,q_t^*) = 0$ was initially considered by ~\cite{Bismut1972,Bismut1978}, see also \cite{Bensoussan1983Stochastic}, and correspond to a wide range of stochastic optimal control problems. For solving the Hamiltonian system \eqref{eq:Hamil_Case_Hu}, we reformulate it by the following new control problem
\begin{equation}\label{eq:new_sc_contr_Hu}
\inf_{\tilde{p}_0,\{\tilde{q}_t\}_{0\leq t\leq T},\{\tilde{u}_t\}_{0\leq t\leq T}}\mathbb{E}\Big[|-h_x(\tilde{x}_T)-\tilde{p}_T|^2 + \displaystyle \lambda
\int_{0}^{T} H_u(t,\tilde{x}_t,\tilde{u}_t,\tilde{p}_t,\tilde{q}_t)^2 \mathrm{d} t \Big],
\end{equation}
\begin{equation*}
\begin{array}{l}
\mbox{s.t. } \tilde{x}_t = x_0 + \displaystyle \int_{0}^{t} b(t,\tilde{x}_s,\tilde{u}_s)\mathrm{d} s + \int_{0}^{t} \sigma(t,\tilde{x}_s,\tilde{u}_s)\mathrm{d} W_s, \vspace{1ex} \\
\hspace{1.8em} \tilde{p}_t = \tilde{p}_0 - \displaystyle \int_{0}^{t} H_x(s,\tilde{x}_s,\tilde{u}_s,\tilde{p}_s,\tilde{q}_s) \mathrm{d} s + \int_{0}^{t} \tilde{q}_s\mathrm{d} W_s,
\end{array}
\end{equation*}
where $\tilde{p}_0,\{\tilde{q}_t\}_{0\leq t\leq T},\{\tilde{u}_t\}_{0\leq t\leq T}$ are the controls and $\lambda$ is the hyper-parameter. As long as the cost functional \eqref{eq:new_sc_contr_Hu} converges to 0, the 4-tuple $(\tilde{x}_t, \tilde{u}_t, \tilde{p}_t, \tilde{q}_t)$ converges to $(x_t^*,u_t^*, p_t^*,q_t^*)$.

The discrete Euler scheme is given as
\begin{equation}\label{eq:Euler_Hu}
\left\{
\begin{array}{l}
\tilde{x}^{\pi}_{t_{i+1}} = \tilde{x}^{\pi}_{t_i} + b(t_i,\tilde{x}^{\pi}_{t_i},\tilde{u}^{\pi}_{t_i})\Delta t_i + \sigma(t_i,\tilde{x}^{\pi}_{t_i},\tilde{u}^{\pi}_{t_i})\Delta W_{t_i}, \vspace{1ex} \\
\tilde{p}^{\pi}_{t_{i+1}} = \tilde{p}^{\pi}_{t_i}-H_x(t_i,\tilde{x}^{\pi}_{t_i},\tilde{u}^{\pi}_{t_i},\tilde{p}^{\pi}_{t_i},\tilde{q}^{\pi}_{t_i})\Delta t_i + \tilde{q}^{\pi}_{t_i}\Delta W_{t_i}, \vspace{1ex} \\
\tilde{x}^{\pi}_0 = x_0, \qquad \tilde{p}^{\pi}_0=\tilde{p}_0,
\end{array}
\right.
\end{equation}
and the loss function is
\begin{equation}\label{eq:cost_fun_Hu}
\mbox{loss} = \dfrac{1}{M} \sum_{j=1}^{M}\Big[ |-h_x(\tilde{x}^{\pi}_T)-\tilde{p}^{\pi}_T|^2 + \lambda \sum_{i=0}^{N-1} H_u(t, \tilde{x}_{t_i}^{\pi},\tilde{u}_{t_i}^{\pi},\tilde{x}_{p_i}^{\pi},\tilde{q}_{t_i}^{\pi})^2 \Big],
\end{equation}
where the time-divided coefficient $T/N$ is merged into the coefficient $\lambda$.

Different from Algorithms \ref{alg:1}, we regard the two processes $\{\tilde{q}_{t_i}^{\pi}, \tilde{u}_{t_i}^{\pi}\}_{0\leq i\leq N-1}$ as controls, which means that two neural networks (2-NNets) should be constructed to simulate $\tilde{q}_{t_i}^{\pi},\tilde{u}_{t_i}^{\pi}$, respectively. And similar with Algorithms \ref{alg:1}, we regard $\tilde{q}_{t_i}^{\pi},\tilde{u}_{t_i}^{\pi}$ as feedback controls of the state $\tilde{x}_{t_i}^{\pi}$ and the time $t_i$, then we construct a common neural network for all time steps, respectively,
\begin{equation}\label{eq:pro_fun_q_Hu}
\left.
\begin{array}{l}
\tilde{q}^{\pi}_{t_i} = \phi^1(t_i, \tilde{x}^{\pi}_{t_i};\theta^q), \vspace{1ex} \\
\tilde{u}^{\pi}_{t_i} = \phi^2(t_i, \tilde{x}^{\pi}_{t_i};\theta^u).
\end{array}
\right.
\end{equation}
The 2-NNets contain both one $(n+1)-$dim input layers and three $(n+10)-$dim hidden layers, the output layers are $(n\times d)-$dim for $\tilde{q}_{t_i}^{\pi}$  and  $k-$dim for  $\tilde{u}_{t_i}^{\pi}$ respectively. The loss function is given as \eqref{eq:cost_fun_Hu}.

The pseudo-code is given as Algorithm \ref{alg:3}.
\begin{algorithm}[H]
\renewcommand{\thealgorithm}{2}
    \caption{Numerical algorithm with 2-NNets}
    \label{alg:3}
    \begin{algorithmic}[1]
        \Require The Brownian motion $ \Delta W_{t_i} $, initial parameters $ (\theta^{q,0},\theta^{u,0},\tilde{p}_0^{0,\pi}) $, learning rate $ \eta $ and hyper-parameter $\lambda$;
        \Ensure The precesses $ \tilde{x}^{l,\pi}_{t_i} $ and $ \tilde{p}^{l,\pi}_{T} $.
        \For { $ l = 0 $ to $ maxstep $}
        \State $ \tilde{x}_{0}^{l,\pi} = x_{0} $, $ \tilde{p}_{0}^{l,\pi} = \tilde{p}_{0}^{l,\pi} $ and $H=0$
        \For { $ i = 0 $ to $ N-1 $}
        \State $\tilde{q}^{l,\pi}_{t_i} = \phi^1(t_i,\tilde{x}^{l,\pi}_{t_i};\theta^{q,l});$
        \State $\tilde{u}^{l,\pi}_{t_i} = \phi^2(t_i,\tilde{x}^{l,\pi}_{t_i};\theta^{u,l});$
        \State $\tilde{x}^{l,\pi}_{t_{i+1}} = \tilde{x}^{l,\pi}_{t_i} + b(t_i,\tilde{x}^{l,\pi}_{t_i},\tilde{u}^{l,\pi}_{t_i})\Delta t_i + \sigma(t_i,\tilde{x}^{l,\pi}_{t_i},\tilde{u}^{l,\pi}_{t_i})\Delta W_{t_i};$
        \State $\tilde{p}^{l,\pi}_{t_{i+1}} = \tilde{p}^{l,\pi}_{t_i}-H_x(t_i,\tilde{x}^{l,\pi}_{t_i},\tilde{u}^{l,\pi}_{t_i},\tilde{p}^{l,\pi}_{t_i},\tilde{q}^{l,\pi}_{t_i})\Delta t_i + \tilde{q}^{l,\pi}_{t_i}\Delta W_{t_i};$
        \State $H = H + H_u(t_i,\tilde{x}_{t_i}^{l,\pi},\tilde{u}_{t_i}^{l,\pi},\tilde{p}_{t_i}^{l,\pi},\tilde{q}_{t_i}^{l,\pi})^2$
        \EndFor
        \State $ J(\tilde{u}^{l,\pi}(\cdot))=\dfrac{1}{M} \sum_{j=1}^{M}\Big[\dfrac{T}{N} \sum_{i=0}^{N-1}f(t_i,\tilde{x}_{t_i}^{l,\pi},\tilde{u}_{t_i}^{l,\pi}) + h(\tilde{x}_T^{l,\pi}) \Big]; $
        \State $ \mbox{loss} = \dfrac{1}{M} \sum_{j=1}^{M}\Big[ |-h_x(\tilde{x}^{l,\pi}_T)-\tilde{p}^{l,\pi}_T|^2 + \lambda H \Big]; $
        \State $ (\theta^{q,l+1},\theta^{u,l+1},\tilde{p}_0^{l+1,\pi})=(\theta^{q,l}, \theta^{u,l}, \tilde{p}_0^{l,\pi})-\eta\nabla\mbox{loss}. $
        \EndFor
    \end{algorithmic}
\end{algorithm}

In Algorithm \ref{alg:3}, instead of solving the maximum condition explicitly or approximately, we only need to consider the control condition $H_u(t,x,u,p,q)=0$. Thus Algorithm \ref{alg:3} can deal with a wide range of high-dimensional problems more effectively even if the optimal control $\tilde{u}$ can not be solved explicitly.

\subsection{Algorithm 3: Numerical Algorithm with explicit expression of \texorpdfstring{$\bar{H}$}{}}\label{subsec:legendre}

The above mentioned  Algorithm \ref{alg:3} provides a method for solving a general kind of high-dimensional stochastic optimal control problems when the optimal control $\tilde{u}$ has not an explicit solution. In this subsection, we introduce another algorithm for solving high-dimensional cases. We show that as long as the function $\bar{H}$ defined by \eqref{eq:H_bar_definition} is known, we can solve a class of high-dimensional stochastic optimal control problems through the deep-learning method established in our previous work \cite{Peng_FBSDE_numerical}.

In more details,
we consider the stochastic optimal control problem \eqref{eq:sc_object} and assume that $\bar{H}$ defined in \eqref{eq:H_bar_definition} is given.
Then the corresponding Hamiltonian system with $\bar{H}$ is given as follows
\begin{equation}\label{eq:FBSDE_s}
    \left\{
    \begin{array}{l}
    \mathrm{d} x_t^* = \bar{H}_p(t,x^*_t,p_t^*,q_t^*)\mathrm{d} t + \bar{H}_q(t,x^*_t,p_t^*,q_t^*)\mathrm{d} W_t, \vspace{1ex} \\
    \mathrm{d} p_t^* = -\bar{H}_x(t,x^*_t,p_t^*,q_t^*)\mathrm{d} t + q_t^*\mathrm{d} W_t, \qquad t\in[0,T], \vspace{1ex} \\
    x^*_0 = x_0, \qquad p_T^*=-h_x(x^*_T),
    \end{array}
    \right.
\end{equation}
which is essentially an \textit{FBSDE}.
If \eqref{eq:FBSDE_s} satisfies the monotonic conditions \cite{Peng1999Fully,Peng2000Problem},
we can solve it in high dimensions through the deep-learning method proposed in \cite{Peng_FBSDE_numerical} for solving the FBSDEs. But different with that in \cite{Peng_FBSDE_numerical}, a single DNN is constructed for all the time points to improve the computing efficiency, as is in Algorithm 1 and 2. Then the optimal state processes $(x^*(\cdot),p^*(\cdot),q^*(\cdot))$ can be obtained.

Finally we get the optimal control $u^*(\cdot)$ through the maximum condition,
\begin{equation}\label{eq-u-H}
  u^*_t = \underset{u\in U}{\arg\max} H(t,x^*_t,u,p^*_t,q^*_t),
\end{equation}
where $H$ is given in \eqref{eq:Ham_con}.

The pseudo-code is shown in Algorithm \ref{alg:4}.
To solve the extremum problem described at line 11 in Algorithm \ref{alg:4},
we can use the similar methods as that in Algorithm\ref{alg:1}, such as BFGS, SLSQP.
However, an important difference is that we need to calculate the extremum problem only once in the whole algorithm, as the optimal state processes $(x^*(\cdot),p^*(\cdot),q^*(\cdot))$ have already been obtained by solving the Hamiltonian system \eqref{eq:FBSDE_s} with deep learning.

\begin{algorithm}[H]
\renewcommand{\thealgorithm}{3}
    \caption{Numerical algorithm with explicit expression of $\bar{H}$}
    \label{alg:4}
    \begin{algorithmic}[1]
    \Require The Brownian motion $ \Delta W_{t_i} $, initial parameters $ (\theta^0,\tilde{p}_0^{0,\pi}) $, learning rate $ \eta $;
    \Ensure The triple precesses $ (\tilde{x}^{l,\pi}_{t_i},\tilde{p}^{l,\pi}_{t_i},\tilde{q}^{l,\pi}_{t_i})$.
    \For { $ l = 0 $ to $ maxstep $}
    \State $ \tilde{x}_{0}^{l,\pi} = x_{0} $, $ \tilde{p}_{0}^{l,\pi} = \tilde{p}_{0}^{l,\pi}; $
    \For { $ i = 0 $ to $ N-1 $}
    \State $\tilde{q}^{l,\pi}_{t_i} = \phi(t_i, \tilde{x}^{l,\pi}_{t_i},\tilde{p}^{l,\pi}_{t_i};\theta^{l});$
    \State $\tilde{x}^{l,\pi}_{t_{i+1}} = \tilde{x}^{l,\pi}_{t_i} + \bar{H}_p(t_i,\tilde{x}^{l,\pi}_{t_i},\tilde{p}^{l,\pi}_{t_i},\tilde{q}^{l,\pi}_{t_i})\Delta t_i + \bar{H}_q(t_i,\tilde{x}^{l,\pi}_{t_i},\tilde{p}^{l,\pi}_{t_i},\tilde{q}^{l,\pi}_{t_i})\Delta W_{t_i};$
    \State $\tilde{p}^{l,\pi}_{t_{i+1}} = \tilde{p}^{l,\pi}_{t_i}-\bar{H}_x(t_i,\tilde{x}^{l,\pi}_{t_i},\tilde{p}^{l,\pi}_{t_i},\tilde{q}^{l,\pi}_{t_i})\Delta t_i + \tilde{q}^{l,\pi}_{t_i}\Delta W_{t_i};$
    \EndFor
    \State $ \mbox{loss} = \dfrac{1}{M} \sum_{j=1}^{M}\Big[ |-h_x(\tilde{x}^{l,\pi}_T)-\tilde{p}^{l,\pi}_T|^2 \Big]; $
    \State $ (\theta^{l+1},\tilde{p}_0^{l+1,\pi})=(\theta^{l},\tilde{p}_0^{l,\pi})-\eta\nabla\mbox{loss}; $
    \EndFor
    \State $ \tilde{u}^{l,\pi}_{t_i} = \arg\max_{u\in U}H(t_i,\tilde{x}^{l,\pi}_{t_i},u,\tilde{p}^{l,\pi}_{t_i},\tilde{q}^{l,\pi}_{t_i}) $;
    \State $ J(\tilde{u}^{l,\pi}(\cdot))=\dfrac{1}{M} \sum_{j=1}^{M}\Big[\dfrac{T}{N} \sum_{i=0}^{N-1}f(t_i,\tilde{x}_{t_i}^{l,\pi},\tilde{u}_{t_i}^{l,\pi}) + h(\tilde{x}_T^{l,\pi}) \Big]. $
  \end{algorithmic}
\end{algorithm}

The choice of algorithm can be determined according to the nature of the function $H$ in the maximum condition and the character of the stochastic optimal control problem. When $\tilde{u}$ can be solved explicitly, it's better to choose Algorithm 3, then the optimal stochastic control problem is degenerated into the problem of solving FBSDEs. When $\tilde{u}$ does not have an explicit representation, all of the three proposed algorithms can be chosen. If the optimal control problem satisfies the conditions of Algorithm \ref{alg:3} or 3, then Algorithm \ref{alg:3} or 3 will be better alternatives, otherwise Algorithm \ref{alg:1} could be adopted. Besides, our proposed framework and algorithms proposed in Section \ref{sec:numerical_alg} can be extended to deal with stochastic optimal control problems where the control domain is non-convex or the state equations are described by fully coupled FBSDEs \cite{Peng1993Backward,hu2018a}. More details will be discussed in Appendix \ref{appendix:non-convex}.

\section{Numerical results}\label{sec:numerical_results}

In this section, we apply our proposed algorithms in solving some stochastic optimal control problems and show the numerical results. Firstly, we give a low-dimensional example to compare the performance of the three algorithms. Then we show some high dimensional cases and give the corresponding results. If not specially mentioned, we set the dimensions $n=k$. And all the examples in this section are calculated with the number of time-points $N=25$ and the batch size of 64. The learning rates and the sample number in the training set varies to adjust to different examples. The sample number in the test set is set to be $512$. The ReLU activation function and Adam optimizer are adopted in the network architectures and the data is normalized before each layer. In order to get more general results, we calculate the means of numerical results from 10 independent runs of each algorithm.
The numerical experiments are performed in PYTHON on a LENOVO computer with a 2.40 Gigahertz (GHz) Inter Core i7 processor and 8 gigabytes (GB) random-access memory (RAM).

\subsection{A low-dimensional example}\label{subsec:low-dim}

Consider the following LQ control system:
\begin{equation}\label{eq:exmp:2_control}
\left\{
\begin{array}{l}
\mathrm{d} x_t = (-\dfrac{1}{4}x_t + u_t)\mathrm{d} t + (\dfrac{1}{5}x_t + u_t) \mathrm{d}  W_t, \vspace{1ex} \\
x(0) = x_0,
\end{array}
\right.
\end{equation}
and the cost functional is defined as
\begin{equation}\label{eq:linear_cost_func_ex1}
  J(0,x_0;u(\cdot)) = \displaystyle \mathbb{E}\left\{ \dfrac{1}{2}
      \int_{0}^{T} [\left\langle \dfrac{1}{2}x_t, x_t \right\rangle + \left\langle 2u_t, u_t \right\rangle]\mathrm{d} t + \dfrac{1}{2} \left\langle Qx_T, x_T \right\rangle
   \right\},
\end{equation}
where we set the dimensions $n=k, d=1$ and $ Q $ is deterministic matrix taking value in $ \mathbb{R}^{n\times n} $. The control domain is $ U = \mathbb{R}^{n} $, and the Hamiltonian $ H $ is
\begin{equation*}
\begin{matrix}
    H(t,x,u,p,q) = \left\langle p, -\dfrac{1}{4}x + u \right\rangle + \left\langle q,\dfrac{1}{5}x + u \right\rangle - \dfrac{1}{4}\left\langle x,x \right\rangle - \left\langle u,u \right\rangle.
\end{matrix}
\end{equation*}
The explicit representation of the optimal control $u^*$ is
\[
\begin{matrix}
    u^* = \dfrac{1}{2}(p+q).
\end{matrix}
\]
Therefore the function $\bar{H}$ has an explicit form
\begin{equation*}
\begin{matrix}
    \bar{H}(x,p,q) = -\dfrac{1}{2}\left\langle \dfrac{1}{2}x, x \right\rangle + \left\langle p, -\dfrac{1}{4}x \right\rangle
  +\left\langle q, \dfrac{1}{5}x \right\rangle+\dfrac{1}{4}\left\langle p+q, p+q \right\rangle.
\end{matrix}
\end{equation*}

The corresponding Hamiltonian system is
\begin{equation}\label{eq:ex2_H}
\left\{
\begin{array}{l}
\mathrm{d} x_t^* = (-\dfrac{1}{4}x_t^* +u_t^*) \mathrm{d} t + (\dfrac{1}{5}x_t^* +u_t^*) \mathrm{d}  W_t, \vspace{1ex} \\
-\mathrm{d} p_t^* = (-\dfrac{1}{2}x_t^*-\dfrac{1}{4}p_t^* +\dfrac{1}{5}q_t^*)\mathrm{d} t - q_t^* \mathrm{d}  W_{t}, \vspace{1ex} \\
x_0^* = x_0, p_T^* = -Qx_T^*, \vspace{1ex} \\
u_t^* = \dfrac{1}{2}(p_t^*+q_t^*).
\end{array}
\right.
\end{equation}
We set $ x_0=1 $ in this example. It can be verified that equation \eqref{eq:ex2_H} satisfies the monotonic condition and has a unique solution $ (x^*(\cdot),p^*(\cdot),q^*(\cdot)) $. The reader can check it in Appendix \ref{appendix:Existence FBSDEs}.

Supposing the solution of FBSDE \eqref{eq:ex2_H} is in the following form:
\begin{align*}
  p_t^* = -K_tx_t^*, \qquad q_t^* = -M_tx_t^*.
\end{align*}
Combing it with \eqref{eq:ex2_H}, we obtain a Riccati equation
\begin{equation}\label{eq:exmp:2_Riccati}
  \begin{cases}
    \dot{K}_t-\dfrac{1}{2}K_t^2-\dfrac{1}{2}K_t+(\dfrac{1}{5}E_n-\dfrac{1}{2}K_t)M_t+\dfrac{1}{2}E_n=0,\vspace{1ex}\\
    \dfrac{1}{2}K_t^2-\dfrac{1}{5}K_t+\dfrac{1}{2}K_tM_t+M_t=0,\vspace{1ex}\\
    K_T = Q,
  \end{cases}
\end{equation}
where $ \dot{K}_t $ is the derivative of $ K_t $ with respect to $ t $, $ E_n $ is an $ n $-order unit matrix. Equation \eqref{eq:exmp:2_Riccati} is a deterministic one-order ordinary differential equation and we can get its numerical solution with the four-order Runge-Kutta methods by using the ODE45 method in Matlab (ODE45 in brief for easy expression). Therefore the numerical solutions of ODE45 is used as a benchmark to be compared with that of our algorithms in the LQ control problem.

Firstly, we set $ T=0.1 $, $ Q=E_n $. The numerical solution of equation \eqref{eq:exmp:2_Riccati} with ODE45 is $ K_0 = 0.9586E_n $ for each dimension, and the value of $ p_0^* $ with ODE45 is
\begin{align*}
  p_0^* = -K(0)x_0 = -0.9586,
\end{align*}
i.e. $ p_0^* $ is a $ n $-dim vector with all its elements equal to $ -0.9586 $.

We give the numerical results comparison of our three proposed algorithms in Table \ref{tab:comparison_u_LQ} with 2000 iteration steps when $ n=5$. And the numerical solution $ p_0^* = -0.9586 $ with ODE45 is used as the benchmark for calculating the relative errors. In order to measure the performance, the approximate solution of the optimal control $\tilde{u}$ in Algorithm \ref{alg:1} are calculated through the L-BFGS method though $\tilde{u}$ can be solved explicitly in this example.
\begin{table}[H]
  \centering
  \caption{Comparison of our different algorithms for $ n=5 $}
  \label{tab:comparison_u_LQ}
  \begin{tabular}{|c|c|c|c|c|c|}
    \hline
    Method &  $p_0$&  Cost& Time(s)& Iteration step &Relative error\\
    \hline
    Alg 1 &-0.95734&2.4008& 5867.3 & 2000&0.13\%\\
    \hline
    Alg 2 &-0.95775&2.3900&110.0    & 2000&0.09\%\\
    \hline
    Alg 3 &-0.95863&2.3880&66.5    & 2000&0.003\%\\
    \hline
  \end{tabular}
\end{table}

From the results in Table \ref{tab:comparison_u_LQ}, we can see that Algorithm 3 is the most effective among the three proposed algorithms with a relative error of $0.003\%$ and running time of $66.5$ s.  And when $\tilde{u}$ has an explicit solution, the function $\bar{H}$ can be obtained, then Algorithm 1 will degenerate to a special case of Algorithm 3. Therefore when $\tilde{u}$ can be solved explicitly, Algorithm 3 is the best choice among all the three algorithms. We also see that even if the optimal control $\tilde{u}$  has not an explicit representation, we can still calculate the stochastic optimal control problem with our algorithms with a relative error of less than $0.2\%$. However, the disadvantage of Algorithm \ref{alg:1} is that it is hard to be applied to high dimensional cases when the solution of $\tilde{u}$ is not explicit, as it is very time consuming to get the approximate solution of $\tilde{u}$ with L-BFGS in each iteration step. And from both the relative errors and the running time, Algorithms \ref{alg:3} and 3 demonstrate much more effective performance on comparing with Algorithms \ref{alg:1}. Therefore, when the optimal control $\tilde{u}$ is not explicit, Algorithms \ref{alg:3} or 3 will be a better choice for high dimensional cases when the conditions mentioned in Section \ref{subsec:numerical_solution 2-NNets} or Section  \ref{subsec:legendre} are satisfied. Otherwise Algorithms \ref{alg:1} could be chosen but it is not suitable for high dimensional cases.

\begin{rem}\label{rem:5}
    For the LQ stochastic optimal problem, the existence of optimal control does not require such strong assumptions as Assumptions 1, 2, the readers can refer to assumption L1 on page 301 in \cite{Yong_stochastic_control}. In this situation, the assumptions in the corresponding Theorem 2 will also be as weak as assumption L1, the proof is similar to Corollary \ref{cor:equivalent_FBSDE} and we omit it.
\end{rem}

\subsection{Some high-dimensional stochastic optimal control problems}

In this subsection, we show some numerical results for high-dimensional cases. As Algorithm \ref{alg:1} is not suitable for high-dimensional problems when the optimal control $\tilde{u}$ is not explicit, we mainly show the results of Algorithm \ref{alg:3} with $H_u(t,x,u,p,q) = 0$ and that of Algorithm 3 when we have an explicit expression of $\bar{H}$.

\subsubsection{A high-dimensional LQ stochastic optimal control problem}

\begin{figure}[H]
  \centering
  \includegraphics[scale=0.35]{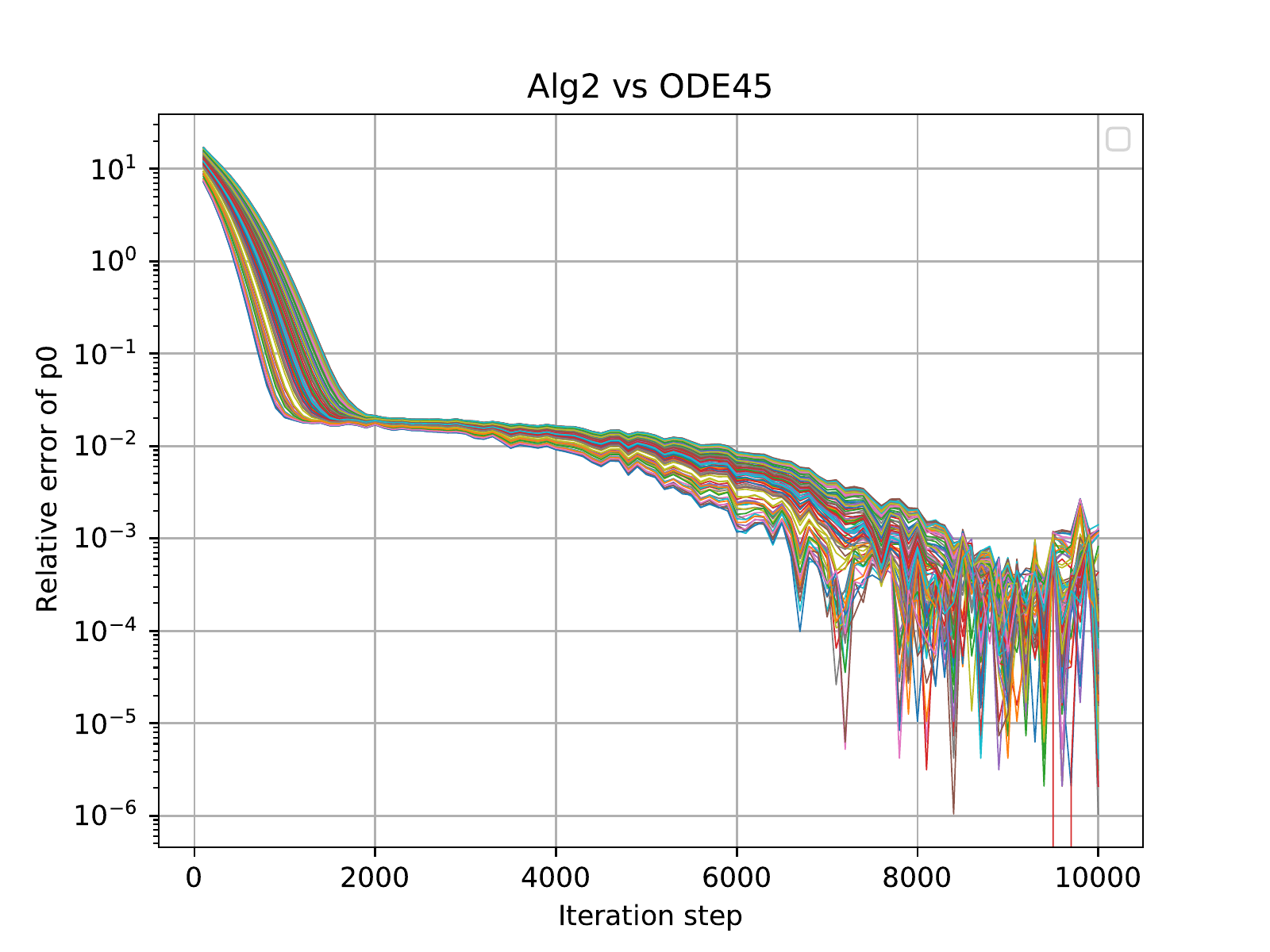}
  \includegraphics[scale=0.35]{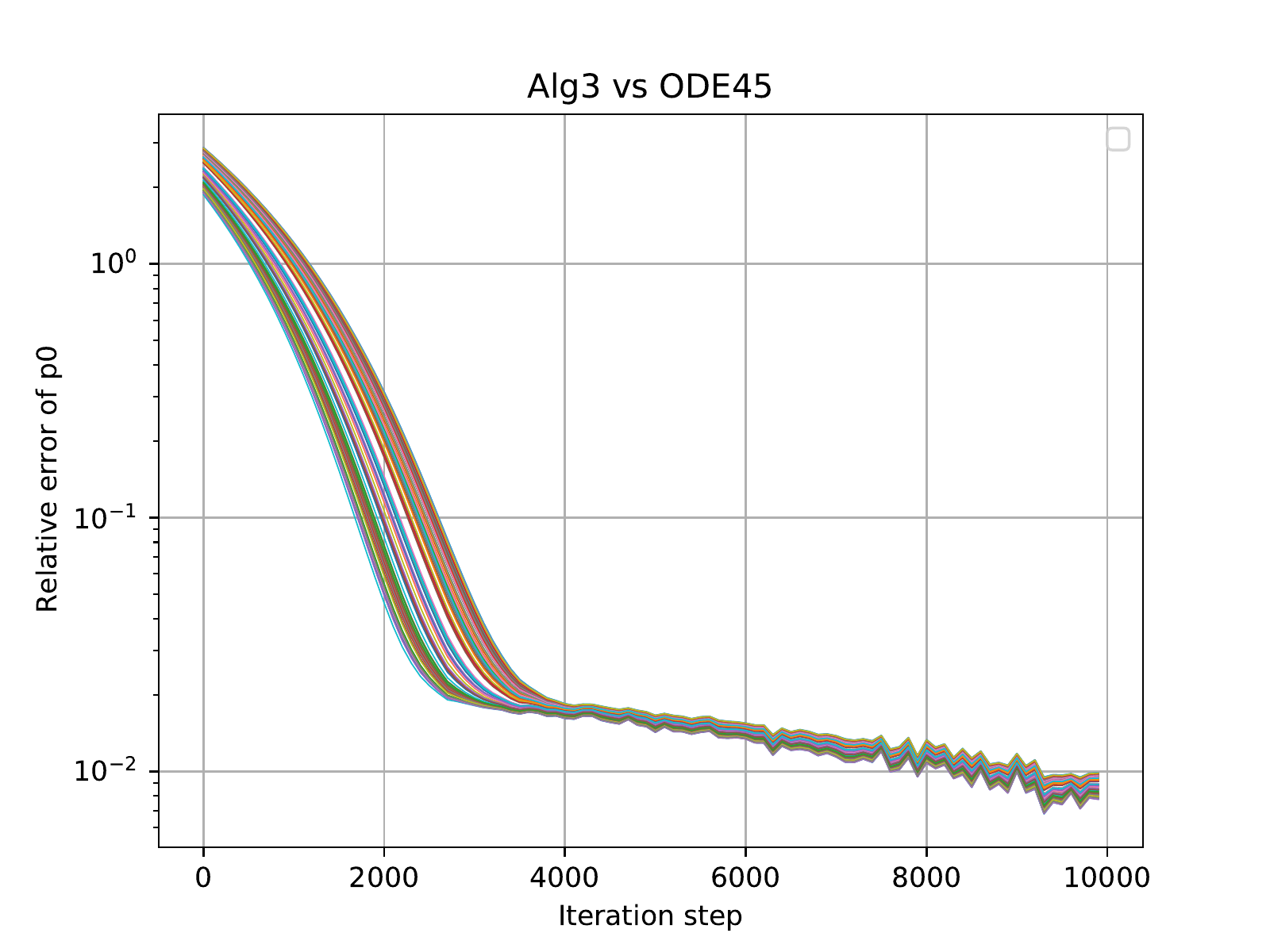}
  \caption{Case $n=100$. The left figure represents the relative errors for each dimension of $ p_0$ with Algorithm \ref{alg:3} on comparing with the results of ODE45, and the right figure shows the relative errors for Algorithm 3.}
  \label{fig:ex_02_d100}
\end{figure}

We first compute the LQ problem mentioned in subsection \ref{subsec:low-dim} for a $n=100$ case. Let $x_0=1.0, T=0.1$. Figure \ref{fig:ex_02_d100} shows the relative errors with different numbers of iteration steps for Algorithm \ref{alg:3} and 3, and the results of ODE45 are used as the benchmark for calculating the relative errors.

\begin{figure}[H]
  \centering
  \includegraphics[scale=0.35]{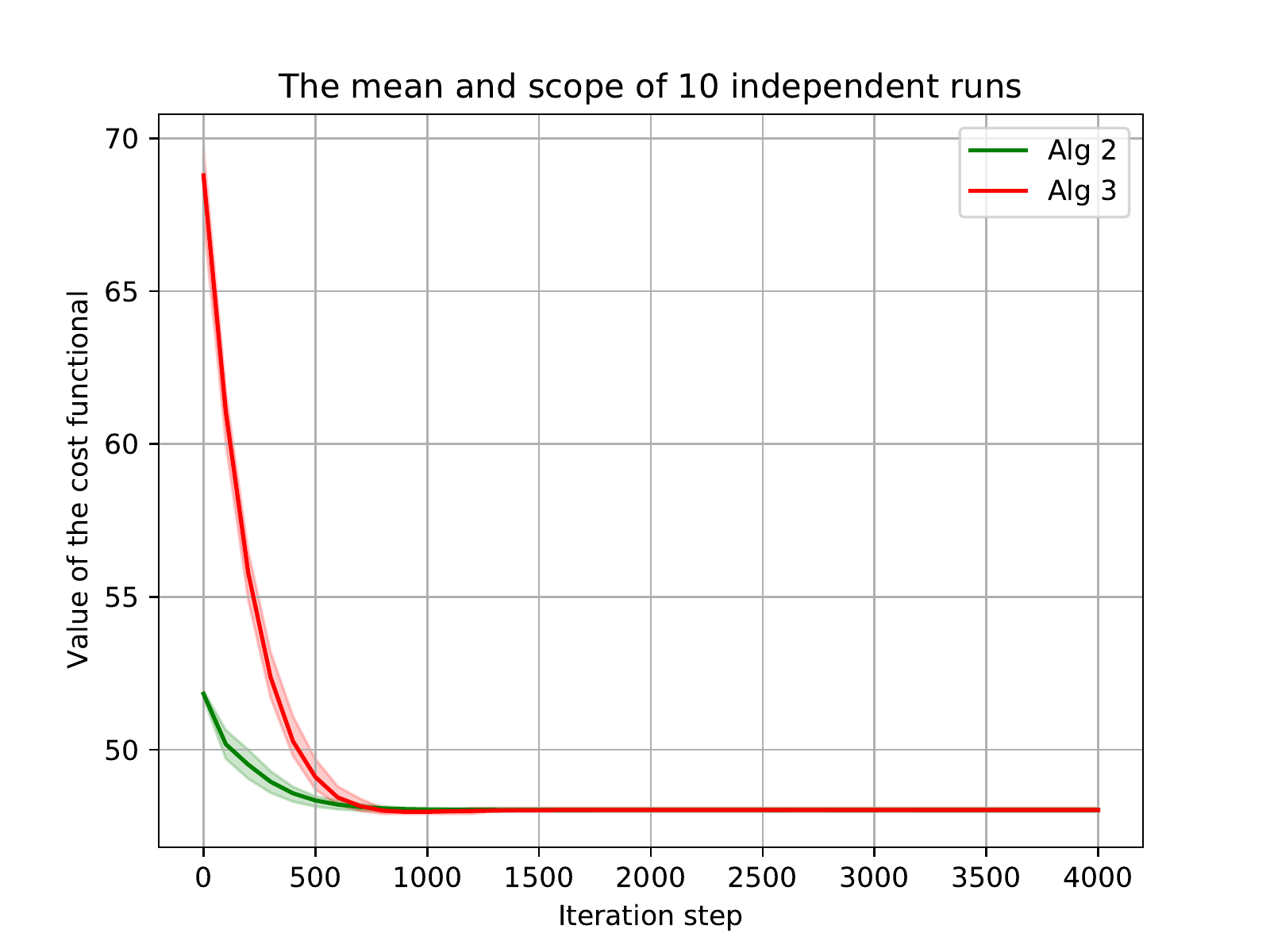}
  \caption{Case $n=100$ and $\lambda=0.05$, the figure shows the mean and scope of the cost functional among 10 independent runs. The green line and scope represent the results of Algorithm \ref{alg:3} with the constraint $H_u(t,x,u,p,q) = 0$. The red line and scope represent the value of cost functional for Algorithm 3 when the optimal control $\tilde{u}$ has an explicit solution. We can see that after 4000 iterations, the values of the cost functionals of Algorithm \ref{alg:3} and 3 are very close and converge to $ 48.056 $.}
\label{fig:ex_02_cost_d100}
\end{figure}

We perform 10 independent runs and the means of the cost functional converges to $ 48.056 $. The curve of the cost functional is shown in Figure \ref{fig:ex_02_cost_d100}.

Now we consider a general case with the terminal $ p_T^* = -\mathbbm{1}_n x_T^* $, where $ \mathbbm{1}_n $ is an $ n $-order matrix with all its elements equal to 1. In this case, the solution $ K_t $ of equation \eqref{eq:exmp:2_Riccati} at time 0 changes with the change of the dimension $ n $.

\begin{figure}[H]
\centering
\includegraphics[width=6cm]{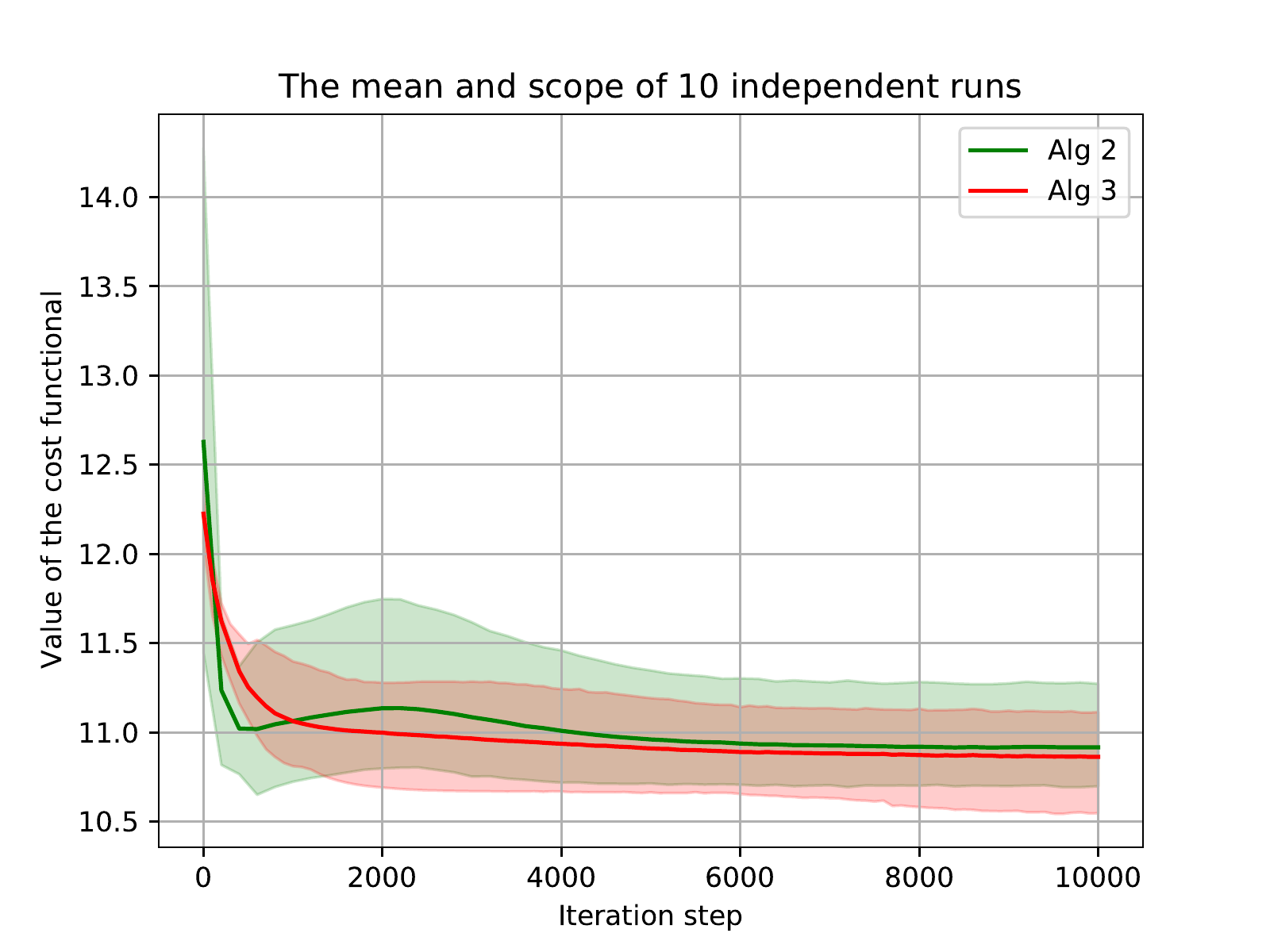}
\caption{Case $n=5$ and $\lambda=0.01$. The figure shows the mean and scope of the cost functional for 10 independent runs. We can see that after 1000 iterations, the value of the cost functional for Algorithm \ref{alg:3} is close to that of Algorithm 3. After 10000 iterations, the mean of the cost functional $J(0, x_0; \tilde{u}(\cdot))$ for Algorithm 2 is $ 10.863 $.}
\label{fig:ex_02_cost_d100_1}
\end{figure}

Similarly, we first calculate the numerical solutions of $p_0$ with ODE45 and take the results as the benchmark. Then we compare our neural network solutions of $p_0$ with that of ODE45. Table \ref{tab:VsMATLAB} shows the comparison results for different dimensions and Figure~\ref{fig:ex_02_cost_d100_1} shows the curve of the cost functional for $ n= 5 $.
We notice that for $n=20$, the results of Alg 3 are not very stable in different runs. In order to get more stable results, we use different neural network for different time points for $n=20$ in Alg 3, but the cost is that more network parameters should be used which need more computation time.

\begin{table}[H]
  \centering
  \caption{Comparison between ODE45 and our algorithms for different dimensions}
  \label{tab:VsMATLAB}
  \resizebox{\textwidth}{20mm}{
  \begin{tabular}{|c|c|c|c|c|c|c|}
    \hline
        \multicolumn{2}{|c|}{} & n=1&  n=2&  n=5&  n=10& n=20 \\
    \hline
    \multicolumn{2}{|c|}{Solution with ODE45 }& -0.9586& -1.8275& -4.3638& -8.5306& -16.821\\
    \hline
    \multirow{2}{*}{Solution with Neural Network}&Alg 2 & -0.9519& -1.8239& -4.3535& -8.4974 & -16.730 \\
    \cline{2-7}
        &Alg 3 & -0.9585 & -1.8276 & -4.3571 & -8.4782& -16.663 \\
        \hline
    \multirow{2}{*}{Relative error}&Alg 2 &0.705\%&    0.197\%& 0.234\%&0.390\% & 0.541\% \\
        \cline{2-7}
        &Alg 3 & 0.007\%& 0.003\%& 0.153\%& 0.615\% & 0.939\%  \\
    \hline
        \multirow{2}{*}{Cost functional}&Alg 2 &0.4812&  1.8375 &10.863&42.850 & 169.80 \\
        \cline{2-7}
        &Alg 3 & 0.4892 & 1.8262 & 10.916 & 43.284 &167.64 \\
    \hline
  \end{tabular}}
\end{table}

\subsubsection{A 100-dim nonlinear control problem}

In this subsection, we compute an example mentioned in \cite{pham2018deep_1,pham2018deep_2}. Consider the following control system
\begin{equation}\label{eq:ex_non_linear_control}
\left\{
\begin{array}{l}
\mathrm{d} x_t = 2u_t\mathrm{d} t + \sqrt{2} \mathrm{d}  W_t, \vspace{1ex} \\
x(0) = x_0
\end{array}
\right.
\end{equation}
and the corresponding cost functional
\begin{equation}\label{eq:non_linear_cost_func_ex1}
\begin{matrix}
    J(0,x_0;u(\cdot)) = \mathbb{E}\left\lbrace \displaystyle \int_{0}^{T} [\left\langle u_t, u_t \right\rangle]\mathrm{d} t + h(x_T)\right\rbrace,
\end{matrix}
\end{equation}
where $ h(x) = \ln(\frac{1}{2}(1+|x|^2)), |x|^2 = \left\langle x, x \right\rangle $. $ W_t $ is a $ n $-dimensional Brownian motion, and the process $ x_t $ is valued in $\mathbb{R}^{n}$. The control domain is $\mathbb{R}$ and the Hamiltonian $ H $ is
\begin{equation*}
  H(t,x,u,p,q) = \left\langle p, 2u \right\rangle + \mbox{tr}(\sqrt{2}q) - \left\langle u,u \right\rangle,
\end{equation*}
the optimal control $ u^* $ can be solved with
\begin{equation*}
  u^* = p,
\end{equation*}
Then we have
\begin{equation*}
  \bar{H}(t,x,p,q) = \left\langle p, 2p \right\rangle + \mbox{tr}(\sqrt{2}q) - \left\langle p,p \right\rangle.
\end{equation*}

\begin{figure}[H]
  \centering
  \includegraphics[scale=0.35]{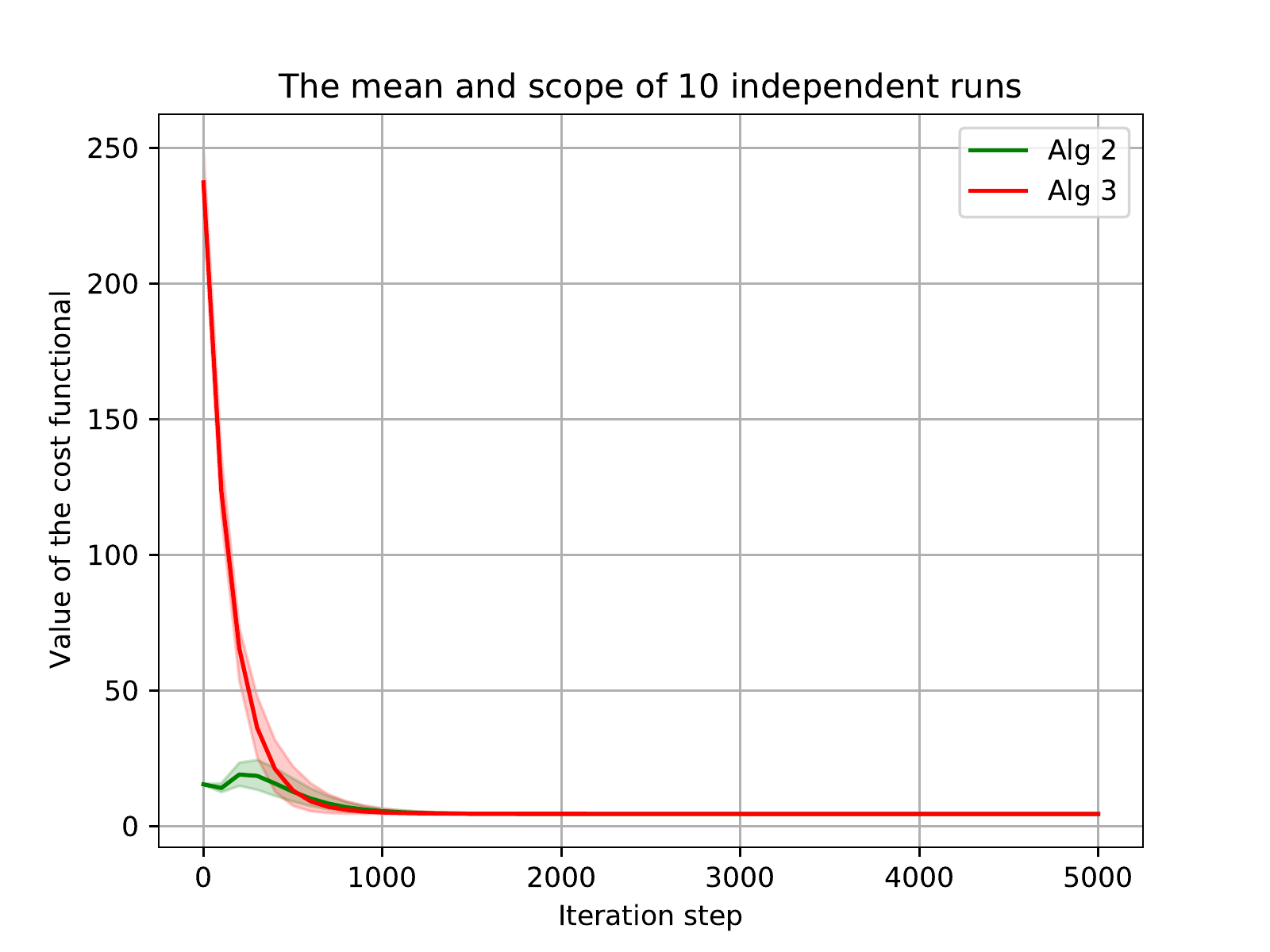}
    \includegraphics[scale=0.35]{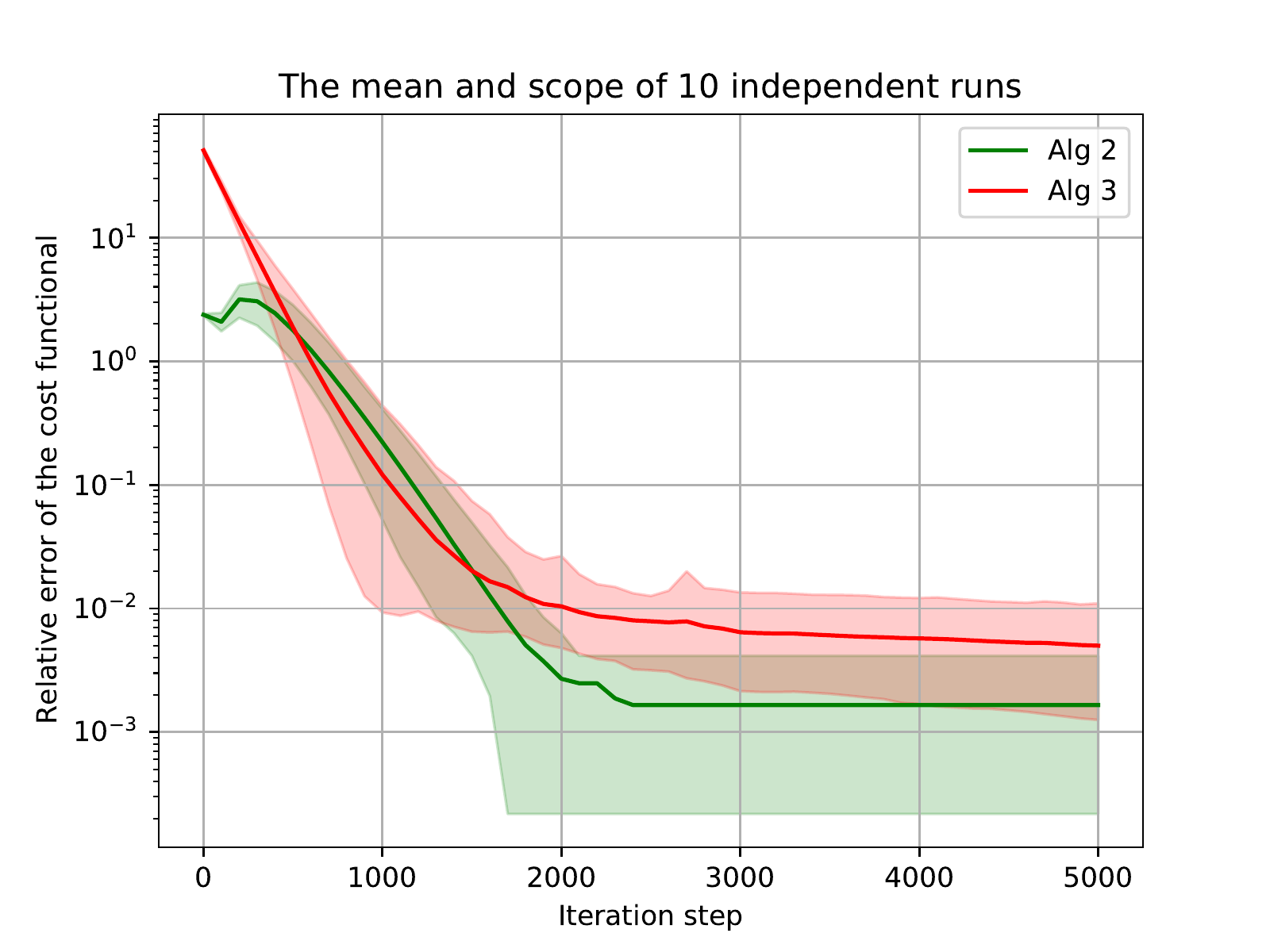}
  \caption{Case $n=100$ and $\lambda=0.05$. The left figure shows the cost functional curve and the right figure shows the mean and scope of relative errors of the cost functional $ J(0, x_0; \tilde{u}(\cdot)) $ among 10 independent runs. We can see that after 5000 iteration steps, the values of the cost functional are very close between Algorithm \ref{alg:3} and 3, and Algorithm \ref{alg:3} shows smaller mean of relative errors.}
  \label{fig:ex_pham_d100}
\end{figure}

The corresponding Hamiltonian system is given as
\begin{equation}\label{eq:non_linear_control}
\left\{
\begin{array}{l}
\mathrm{d} x_t^* = 2u_t^*\mathrm{d} t + \sqrt{2} \mathrm{d}  W_t, \vspace{1ex} \\
-\mathrm{d} p_t^* = -q_t^* \mathrm{d}  W_t, \vspace{1ex} \\
x_0^* = x_0, \qquad p_T^* = -h_{x}(x_T^*), \vspace{1ex} \\
u_t^* = p_t^*
\end{array}
\right.
\end{equation}
where $ p^*_{\cdot} $ and $ q^*_{\cdot} $ are valued in $ \mathbb{R}^{n} $ and $ \mathbb{R}^{n\times n} $ respectively and $n=d=k$. It is known that the cost functional of \eqref{eq:non_linear_cost_func_ex1} can be obtained via the Hopf-Cole transformation, as shown in the following:
\begin{equation*}
  J(0,x_0;u^*(\cdot)) = -\ln(\mathbb{E}[exp(-h(x_0+\sqrt{2}W_T))]), \qquad x_0\in\mathbb{R}^{n}.
\end{equation*}

which can be approximated by Monte Carlo simulation. We set $ x_0 = 0$ and $T=1.0 $ with the dimension $ n=100 $ and the approximate cost functional is 4.591 by performing the Monte Carlo simulation. Figure \ref{fig:ex_pham_d100} shows the values of the cost functional and the relative errors of the cost functional on comparing with the results of the Monte Carlo simulation.

\subsubsection{An application in calculating sub-linear expectation and nonlinear PDE}

Our algorithm can also be used to solve the sub-linear expectation~\cite{peng2019nonlinear}, which is an expansion of the classical expectation theory and can be applied in the situations of model uncertainty. Moreover, the sub-linear expectation is connected to one kind of fully nonlinear PDEs, i.e. the proposed algorithms are also suitable for solving some fully non-linear PDEs.

In this subsection, we calculate the \texorpdfstring{$ G $}{}-expectation which is a special and important sub-linear expectation. The preliminary knowledge of the \texorpdfstring{$ G $}{}-expectation can be seen in ~\cite{peng2019nonlinear}.

For a given function $\varphi$ and a $n$-dimensional $G$-normal distributed variable $X$, we compute the sub-linear expectation $\hat{\mathbb{E}}[\varphi(X)]$. According to the representation of the sub-linear expectation, we have
\begin{equation*}
  \hat{\mathbb{E}}[\varphi(X)] = \sup_{\theta\in\Theta}\mathbb{E}_{\theta}[\varphi(X)] = -\inf_{\theta\in\Theta}\mathbb{E}_{\theta}[-\varphi(X)],
\end{equation*}
which is equivalent to the following stochastic optimal control problem
\begin{equation}\label{eq:G_control}
\left\{
\begin{array}{l}
\mathrm{d} X_t = \theta_t\mathrm{d} W_t, \vspace{1ex} \\
\underline{\sigma}\leq \theta_t\leq\bar{\sigma}, \vspace{1ex} \\
X(0)=0,
\end{array}
\right.
\end{equation}
with the cost functional
\begin{equation*}
  J(\theta(\cdot)) = \inf_{\theta\in\Theta}\mathbb{E}_{\theta}[-\varphi(X)],
\end{equation*}
where $\theta$ is the control and $\Theta$ is the control domain and $ W_t $ is a $ n $-dimensional Brownian motion. \eqref{eq:G_control} is a stochastic optimal control problem with constraint and its corresponding Hamiltonian $ H $ can be easily got by
\begin{equation*}
H(t,x,\theta,p,q)=\theta q
\end{equation*}
with the constraint $ \underline{\sigma}\leq \theta\leq\bar{\sigma} $. The optimal control $ \theta^* $ is
\begin{equation*}
\theta^* = \bar{\sigma}\mathbbm{1}_{q\geq 0} + \underline{\sigma}\mathbbm{1}_{q<0},
\end{equation*}
and $\bar{H}$ is given as
\begin{equation*}
\bar{H}(t,x,p,q)=(\bar{\sigma}\mathbbm{1}_{q\geq 0} + \underline{\sigma}\mathbbm{1}_{q<0}) q.
\end{equation*}

The Hamiltonian system can be got as follows
\begin{equation}\label{eq:G_Hamilton}
\left\{
\begin{array}{l}
\mathrm{d} X_t^* = \theta_t^*\mathrm{d} W_t, \vspace{1ex} \\
-\mathrm{d} p_t^* = -q_t^*\mathrm{d} W_t, \vspace{1ex} \\
X_0^*=0, \qquad p_T^*=\varphi_{x}(X_T^*), \vspace{1ex} \\
\theta_t^*=\bar{\sigma}\mathbbm{1}_{q_t^*\geq 0} + \underline{\sigma}\mathbbm{1}_{q_t^*<0},
\end{array}
\right.
\end{equation}
which connects to one fully nonlinear PDE
\begin{equation}\label{Gheat}
  \partial_t u_t - G(\partial_{xx}^2 u_t) = 0, \qquad u|_{t=0} = \varphi,
\end{equation}
where
\begin{equation*}
    \left\{
    \begin{array}{l}
        u(t,x) = \hat{\mathbb{E}} [\varphi(x + \sqrt{t} X)], \vspace{1ex} \\
        G(A) = \dfrac{1}{2} \hat{\mathbb{E}} [\langle AX, X \rangle].
    \end{array}
    \right.
\end{equation*}
In particular, $\hat{\mathbb{E}} [\varphi(X)] = u(1, 0)$. The parabolic PDE \eqref{Gheat} is called a $G$-heat equation. Therefore we could solve one kind of fully nonlinear PDEs.

In problem \eqref{eq:G_Hamilton}, for the cases where the function $\varphi$ is convex, we can solve it through our framework even though the assumptions in Corollary \ref{cor:equivalent_FBSDE} are not satisfied. The proof is similar to Corollary \ref{cor:equivalent_FBSDE} and we omit it.

\begin{figure}[H]
  \centering
  \includegraphics[scale=0.35]{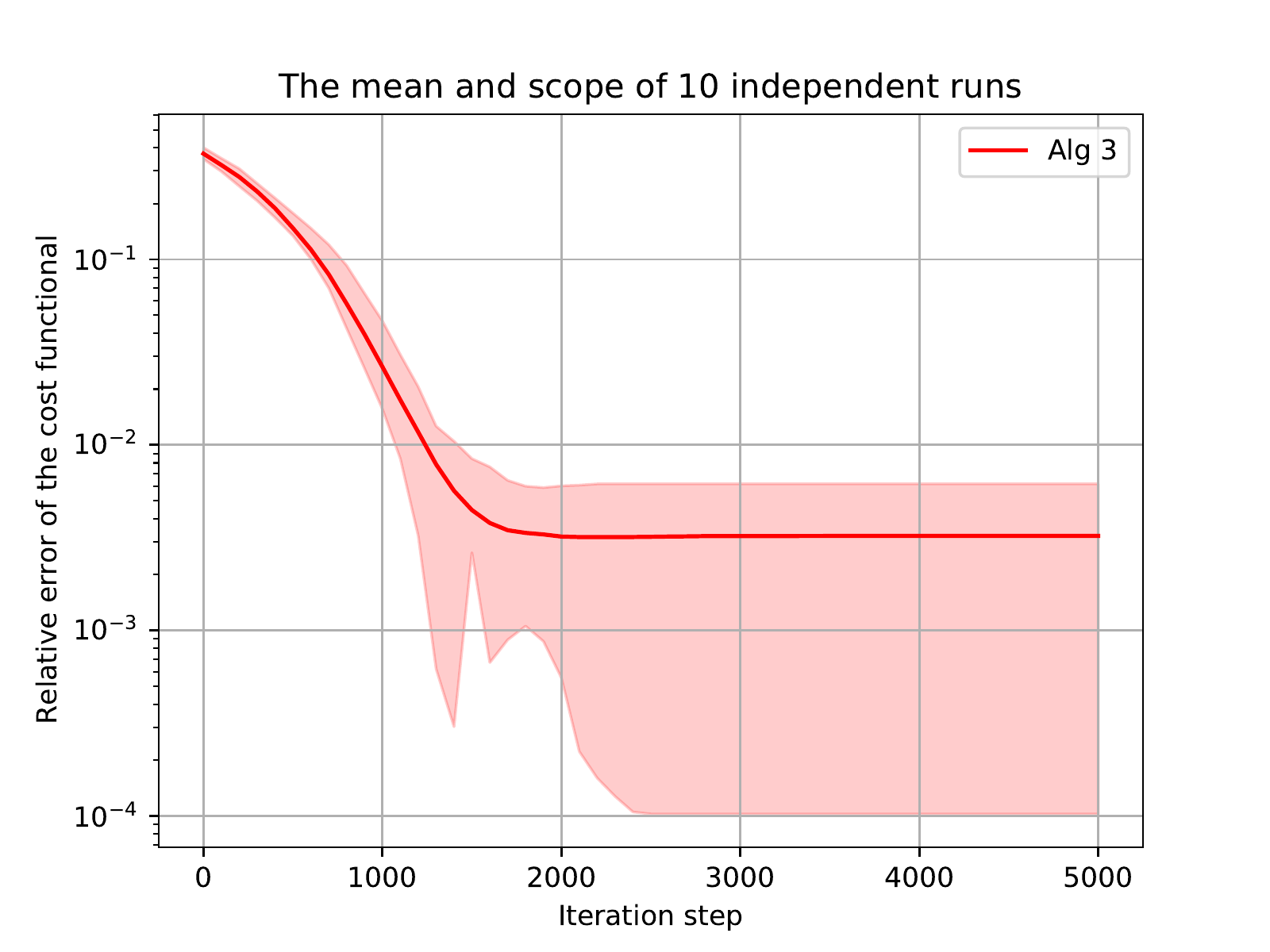}
  \caption{Case $n=100$. This figure shows the mean and scope of relative errors of the cost functional $ J(0,x_0;\theta^*(\cdot)) $ among 10 independent runs. After 5000 iterations, the value of $ J(\theta^*(\cdot)) $ converges to 400.268.}
  \label{fig:ex_G_d100}
\end{figure}

Let $ n = 100 $, $ X \overset{d}{=}N(\{0\}\times\Theta)= N({0}\times[\underline{\sigma}^2,\bar{\sigma}^2]^n)  $ with $ \underline{\sigma}=1, \bar{\sigma} = 2 $ and $ \varphi(x) = x^2 $. The exact value of $ \hat{\mathbb{E}}[\varphi(X)] $ is equal to 400 which is taken as the benchmark for measuring the performance of our algorithm. As the partial of the Hamiltonian $H_{\theta}$ doe not contain any term with the control $\theta$, thus Algorithm \ref{alg:3} with the constraint $H_\theta(t,x,\theta,p,q) = 0$ can not be applied to this example, therefore we mainly show the results of Algorithm 3. Figure \ref{fig:ex_G_d100} shows the mean and scope of relative errors of the cost functional among 10 independent runs.

\subsubsection{An example of the control \texorpdfstring{$\tilde{u}$}{} without explicit solution}\label{ssub:no_solution}

\begin{figure}[H]
    \centering
    \includegraphics[scale=0.35]{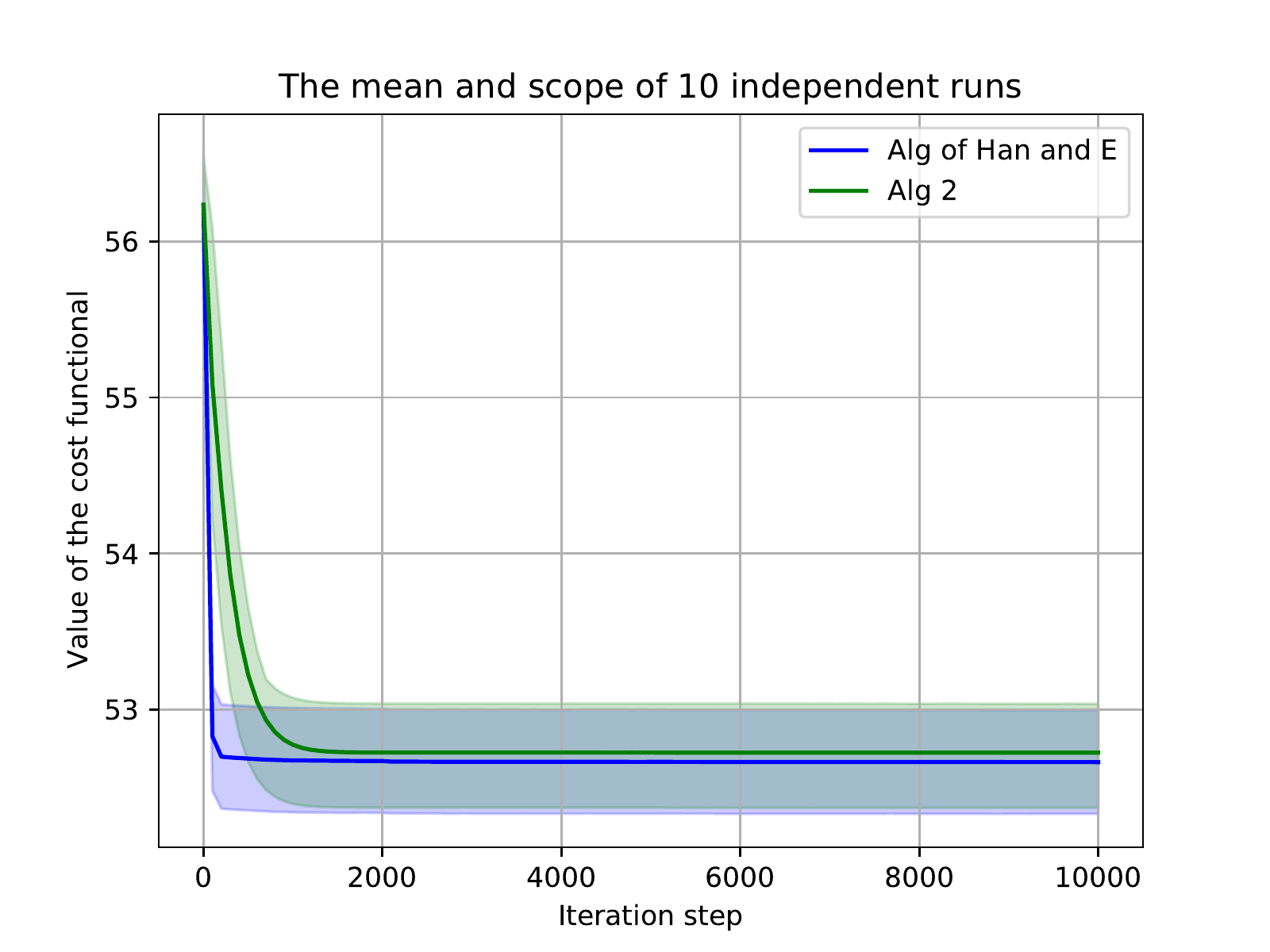}
    \includegraphics[scale=0.35]{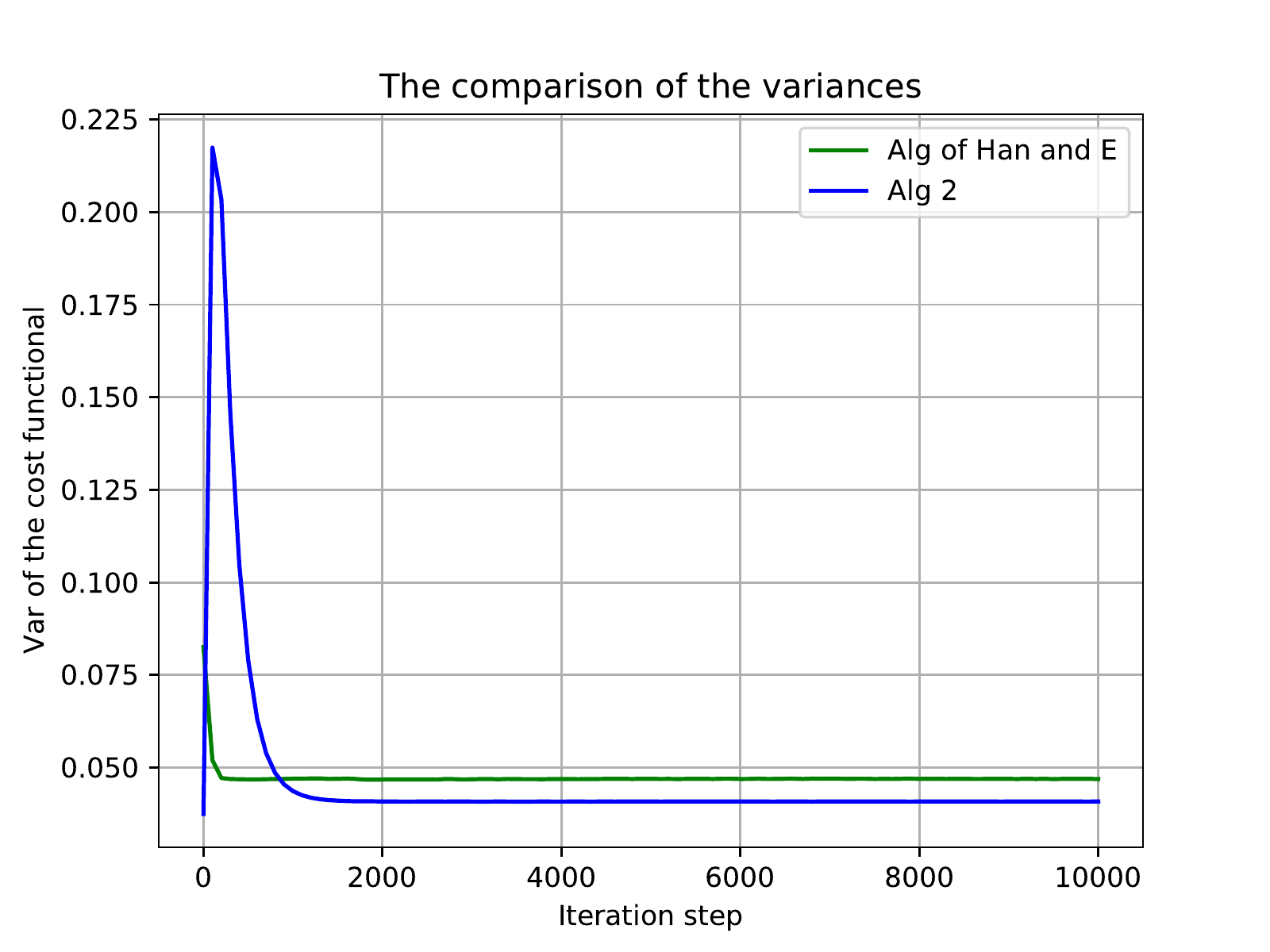}
    \caption{Case $n=100$ and $\lambda=0.1$. The left figure shows the comparison between our algorithm and that of Han and E on the mean and scope of the cost functional among 10 independent runs. The right figure shows the variances comparison. We can see that the value of the cost functional for the two algorithms are very close from the left figure and our algorithm has lower variances from the right figure.}
    \label{fig:ex_05_d100}
\end{figure}

In this subsection, we show an example whose optimal control $\tilde{u}$ doesn't have an explicit form. Consider the following stochastic control problem,
\begin{equation}\label{eq:transcendental}
\left\{
\begin{array}{l}
\mathrm{d} x_t = \sin u_t \mathrm{d} t + x_t \mathrm{d} W_t, \vspace{1ex} \\
x(0) = x_0,
\end{array}
\right.
\end{equation}
with cost functional
\begin{equation*}\label{eq:cost_sec5}
\begin{matrix}
    J(u(\cdot)) = \mathbb{E} \left\lbrace \displaystyle \int_{0}^{T} \langle u_t, u_t \rangle \mathrm{d} t + \dfrac{1}{2} \langle x_T, x_T \rangle \right\rbrace,
\end{matrix}
\end{equation*}
where the control domain is $U=\mathbb{R}^n$, $ W_t $ is a $ n $-dimensional Brownian motion. The Hamiltonian $H$ is

\begin{equation*}
    H(t,x,u,p,q) = \langle p, \sin u \rangle + \langle q, x \rangle - \langle u, u \rangle,
\end{equation*}
which is a multi-dimensional transcendental equation and has not an explicit representation of both the optimal control $\tilde{u}$ and the function $\bar{H}$. The derivative of the Hamiltonian $H$ in $u$ is given as
\begin{equation*}
    H_{u}(t,x,u,p,q) = p \cos u - 2u,
\end{equation*}
which is also a multi-dimensional transcendental equation. The corresponding Hamiltonian system  with the constraint $H_{u}(t,x,u,p,q) = 0 $ is
\begin{equation}\label{eq:transcen_equ_system}
\left\{
\begin{array}{l}
\mathrm{d} x_t^* = \sin u_t^*\mathrm{d} t + x_t \mathrm{d}  W_t, \vspace{1ex} \\
-\mathrm{d} p_t^* = q_t^* \mathrm{d} t -q_t^* \mathrm{d}  W_t, \vspace{1ex} \\
x_0^* = x_0, \qquad p_T^* = -h_{x}(x_T^*), \vspace{1ex} \\
p_t^* \cos u_t^* - 2u_t^* = 0.
\end{array}
\right.
\end{equation}
We set $x_0 = 1.0$, $T = 0.1$ and the hyper-parameter $\lambda = 0.1$. As the optimal control $\tilde{u}$ has not an explicit solution and $\bar{H}$ does not have an explicit expression, thus Algorithm 3 is not suitable for this example. Therefore we mainly give the numerical results of Algorithm \ref{alg:3} in this example. In order to show the convergence of our algorithm, we implement the algorithm of Han and E \cite{han2016deep} which solves the stochastic optimal control problems directly via deep learning, then compare our algorithm with that of Han and E. The comparing results are shown in Figure~\ref{fig:ex_05_d100}. We can see that the value of the cost functional for the two algorithms are very close, but the advantage of our algorithm is that it provides a criterion to decide whether the state-control pair $(\tilde{x}(\cdot), \tilde{u}(\cdot))$ is an optimal pair, that is, whether the value of the loss function equals to 0.

\section{Conclusion}\label{sec:conclusion}
In this paper, we have solved the stochastic optimal control problem from the view of the stochastic maximum principle and proposed three different algorithms via deep learning. We have compared our proposed algorithms through numerical results and pointed out their applicative situations. The numerical results for different examples demonstrate the effectiveness of our proposed algorithms.

\appendix
For readers' convenience, the programming codes of our proposed algorithms is available at the following website: https://github.com/mathfinance-sdu/Deep-solver-for-stochastic-optimal-control-with-SMP. Some preliminaries and the algorithm for solving problems with non-convex control domain are listed in the appendix.

   \section{Problems with non-convex control domain}\label{appendix:non-convex}

  The three algorithms proposed in this paper are mainly aimed at the stochastic optimal control problems with first-order adjoint equations which correspond to the cases with convex control domain. In fact, for the cases with non-convex control domain, we can use the similar method. In this situation, we need to introduce the following second order adjoint equation when $\sigma$ contains the control $u$:
  \begin{equation}\label{eq:2-order-adjoint}
  \left\{
  \begin{array}{l}
  \mathrm{d} P_t^* = -\Big\{b_x(t,x^*_t,u^*_t)\transpose P_t^* + P_t^*b_x(t,x^*_t,u^*_t) + \sum_{j=1}^{d}\sigma_x^j(t,x^*_t,u^*_t)\transpose P_t^*\sigma_x^j(t,x^*_t,u^*_t) \vspace{1ex} \\
  \hspace{5em}+\sum_{j=1}^{d}\Big\{\sigma_x^j(t,x^*_t,u^*_t)\transpose Q_{jt}^* + Q_{jt}^*\sigma_x^j(t,x^*_t,u^*_t)\Big\} \vspace{1ex} \\
  \hspace{5em}+H_{xx}(t,x^*_t,u^*_t,p_t^*,q_t^*)\Big\}\mathrm{d} t + \sum_{j=1}^{d}Q_{jt}^*\mathrm{d} W^{j}_t, \vspace{1ex} \\
  P_T^*=-h_{xx}(x^*_T),
  \end{array}
  \right.
  \end{equation}
  where the \textit{Hamiltonian} $ H $ is defined by
  \begin{equation}\label{eq:Ham_con_1}
  \begin{array}{l}
  H(t,x,u,p,q) = \left\langle  p,b(t,x,u) \right\rangle + \mbox{tr}[q\transpose\sigma(t,x,u)] - f(t,x,u), \vspace{1ex} \\
  \hspace{7em}((t,x,u,p,q))\in[0,T]\times\mathbb{R}^n\times U\times \mathbb{R}^n \times \mathbb{R}^{n\times d},
  \end{array}
  \end{equation}
  and $ (p^*(\cdot),q^*(\cdot)) $ is the solution of \eqref{eq:1-order-adjoint}. In equation \eqref{eq:2-order-adjoint}, the solution is a pair of processes $ (P^*(\cdot),Q^*(\cdot))\in L_{\mathcal{F}}^2(0,T;\mathcal{S}^n)\times(L_{\mathcal{F}}^2(0,T;\mathcal{S}^n))^d $ where $ \mathcal{S}^n = \{ A\in\mathbb{R}^{n\times n}|A\transpose=A \} $.

  Equation \eqref{eq:2-order-adjoint} is also a BSDE with matrix-valued $ (P^*(\cdot),Q^*(\cdot)) $. As with \eqref{eq:1-order-adjoint}, there exists a unique adapted solution $ (P^*(\cdot),Q^*(\cdot)) $ to \eqref{eq:2-order-adjoint} under Assumption \ref{assu:1}. We refer to \eqref{eq:1-order-adjoint} (resp. \eqref{eq:2-order-adjoint}) as the \textit{first-order} (resp. \textit{second-order}) \textit{adjoint equations} and to $ p^*(\cdot) $ (resp. $P^*(\cdot)$) as the \textit{first-order} (resp. \textit{second-order}) \textit{adjoint process}. If $ (x^*(\cdot),u^*(\cdot)) $ is an optimal (resp. admissible) pair, and $ (p^*(\cdot),q^*(\cdot)) $ and $ (P^*(\cdot),Q^*(\cdot)) $ are adapted solutions of \eqref{eq:1-order-adjoint} and \eqref{eq:2-order-adjoint}, respectively, then $ (x^*(\cdot), u^*(\cdot), p^*(\cdot), q^*(\cdot), P^*(\cdot), Q^*(\cdot)) $ is called an \textit{optimal 6-tuple} (resp. \textit{admissible 6-tuple}). The following theorem can be found in \cite{Peng90,Yong_stochastic_control}.

  \begin{thm}\label{thm:Pont_MP_1}
    Let Assumption \ref{assu:1} hold. Let $ (x^*(\cdot),u^*(\cdot)) $ be an optimal pair of \eqref{eq:sc_object}. Then there are pairs of processes
    \begin{equation}
    \left\{
    \begin{array}{l}
    (p^*(\cdot),q^*(\cdot)) \in L_{\mathcal{F}}^2(0,T;\mathbb{R}^n)\times(L_{\mathcal{F}}^2(0,T;\mathbb{R}^n))^d, \\
    (P^*(\cdot),Q^*(\cdot)) \in L_{\mathcal{F}}^2(0,T;\mathcal{S}^n)\times(L_{\mathcal{F}}^2(0,T;\mathcal{S}^n))^d,
    \end{array}
    \right.
    \end{equation}
    where
    \begin{equation}
    \left\{
    \begin{array}{l}
    q^*(\cdot) = (q^*_1(\cdot),\cdots,q^*_d(\cdot)), \qquad Q^*(\cdot) = (Q^*_1(\cdot),\cdots,Q^*_d(\cdot)), \\
    q_j^*(\cdot) \in L_{\mathcal{F}}^2(0,T;\mathbb{R}^n), \qquad Q_j^*(\cdot)\in L_{\mathcal{F}}^2(0,T;\mathcal{S}^n), \qquad 1\leq j\leq d,
    \end{array}
    \right.
    \end{equation}
    satisfy the first-order and second-order adjoint equations \eqref{eq:1-order-adjoint} and \eqref{eq:2-order-adjoint}, respectively, such that a generalized Hamiltonian
    \begin{equation}\label{eq:max_cond_1}
    \mathcal{H}(t,x^*_t,u^*_t) = \max_{u\in U}\mathcal{H}(t,x^*_t,u), \qquad \mbox{a.e. }t\in[0,T], \qquad \mathbb{P}\mbox{-a.s.}
    \end{equation}
    where
    \begin{equation}
    \begin{array}{l}
    \mathcal{H}(t,x,u) \\
    \hspace{1.5em}\triangleq H(t,x,u,p_t,q_t)-\dfrac{1}{2}\mbox{tr}[\sigma(t,x^*_t,u^*_t)\transpose P_t\sigma(t,x^*_t,u^*_t)] \\
    \hspace{2em}+\dfrac{1}{2}\mbox{tr}\Big\{[\sigma(t,x,u)-\sigma(t,x^*_t,u^*_t)]\transpose P_t  [\sigma(t,x,u)-\sigma(t,x^*_t,u^*_t)] \Big\}.
    \end{array}
    \end{equation}
  \end{thm}
  The sufficient conditions for the optimality can be found in \cite{Yong_stochastic_control} and the unique optimal control $u^*$ can be got under some strictly convexity assumptions. Then we have the corresponding new stochastic optimal control problem:
  \begin{equation}\label{eq:new_sc_contr_second}
  \inf_{\tilde{p}_0,\tilde{P}_0,\{\tilde{q}_t, \tilde{Q}_t\}_{0\leq t\leq T}}\mathbb{E}\Big[|-h_x(\tilde{x}_T)-\tilde{p}_T|^2 + |-h_{xx}(\tilde{x}_T)-\tilde{P}_T|^2\Big],
  \end{equation}
  \begin{equation*}
  \begin{array}{l}
  \mbox{s.t. } \tilde{x}_t = x_0 + \displaystyle \int_{0}^{t} b(t,\tilde{x}_s,\tilde{u}_s)\mathrm{d} s + \int_{0}^{t} \sigma(t,\tilde{x}_s,\tilde{u}_s)\mathrm{d} W_s, \vspace{1ex} \\
  \hspace{1.8em} \tilde{p}_t = \tilde{p}_0 - \displaystyle \int_{0}^{t} H_x(s,\tilde{x}_s,\tilde{u}_s,\tilde{p}_s,\tilde{q}_s)\mathrm{d} s + \int_{0}^{t} \tilde{q}_s\mathrm{d} W_s, \vspace{1ex} \\
  \hspace{1.5em} \tilde{P}_t = \tilde{P}_0 - \displaystyle \int_{0}^{t} F(t,\tilde{x}_s,\tilde{u}_s,\tilde{p}_s,\tilde{q}_s,\tilde{P}_s,\tilde{Q}_s)\mathrm{d} t + \int_{0}^{t} \tilde{Q}_s\mathrm{d} W_s,
  \end{array}
  \end{equation*}
  where
  \begin{align*}
    F(t,x,u,p,q,P,Q) = &b_x(t,x,u)\transpose P + Pb_x(t,x,u) + \sum_{j=1}^{d}\sigma_x^j(t,x,u)\transpose P\sigma_x^j(t,x,u) \\
  &+\sum_{j=1}^{d}\Big\{\sigma_x^j(t,x,u)\transpose Q_{j} + Q_{j}\sigma_x^j(t,x,u)\Big\}+H_{xx}(t,x,u,p,q),
  \end{align*}
  and the Euler scheme is
  \begin{equation}\label{eq:for_stoch_diff_eq_dis_2}
  \left\{
  \begin{array}{l}
  \tilde{x}^{\pi}_{t_{i+1}} = \tilde{x}^{\pi}_{t_i} + b(t_i,\tilde{x}^{\pi}_{t_i},\tilde{u}^{\pi}_{t_i})\Delta t_i + \sigma(t_i,\tilde{x}^{\pi}_{t_i},\tilde{u}^{\pi}_{t_i})\Delta W_{t_i}, \vspace{1ex} \\
  \tilde{p}^{\pi}_{t_{i+1}} = \tilde{p}^{\pi}_{t_i}-H_x(t_i,\tilde{x}^{\pi}_{t_i},\tilde{u}^{\pi}_{t_i},\tilde{p}^{\pi}_{t_i},\tilde{q}^{\pi}_{t_i})\Delta t_i + \tilde{q}^{\pi}_{t_i}\Delta W_{t_i}, \vspace{1ex} \\
  \tilde{P}^{\pi}_{t_{i+i}} = \tilde{P}^{\pi}_{t_i}- F(t_i,\tilde{x}^{\pi}_{t_i},\tilde{u}^{\pi}_{t_i},\tilde{p}^{\pi}_{t_i},\tilde{q}^{\pi}_{t_i},\tilde{P}^{\pi}_{t_i},\tilde{Q}^{\pi}_{t_i})\Delta t_i + \tilde{Q}^{\pi}_{t_i}\Delta W_{t_i}, \vspace{1ex} \\
  \tilde{x}^{\pi}_0 = x_0, \qquad \tilde{p}^{\pi}_0=\tilde{p}_0, \qquad \tilde{P}^{\pi}_0=\tilde{P}_0, \vspace{1ex} \\
  \mathcal{H}(t_i,\tilde{x}^{\pi}_{t_i},\tilde{u}^{\pi}_{t_i}) = \max_{u\in U}\mathcal{H}(t_i,\tilde{x}^{\pi}_{t_i},u).
  \end{array}
  \right.
  \end{equation}
  For the second-order case, we need to construct two neural networks at the same time, one for simulating $ \tilde{q}^{\pi}_{\cdot} $ and the other for simulating $ \tilde{Q}^{\pi}_{\cdot} $,
  \begin{equation}\label{eq:pro_fun_qQ_2}
  \left\{
  \begin{array}{l}
  \tilde{q}^{\pi}_{t_i} = \phi^1(t_i,\tilde{x}^{\pi}_{t_i};\theta^1), \vspace{1ex} \\
  \tilde{Q}^{\pi}_{t_i} = \phi^2(t_i,\tilde{x}^{\pi}_{t_i};\theta^2).
  \end{array}
  \right.
  \end{equation}
  These two networks have both one $ (1+n) $-dim input layer, and the output layers are $ (n\times d) $-dim and $ (n\times n\times d) $-dim, respectively. All parameters of the two networks are represented as $ \theta $ and the loss function is defined as
  \begin{equation}\label{eq:cost_fun}
  \mbox{loss} = \dfrac{1}{M} \sum_{j=1}^{M}\Big[ |-h_x(\tilde{x}^{\pi}_T)-\tilde{p}^{\pi}_T|^2 + |-h_{xx}(\tilde{x}^{\pi}_T)-\tilde{P}^{\pi}_T|^2 \Big],
  \end{equation}
  where $ M $ is the number of samples. The pseudo-code for the second-order case is given in Algorithm \ref{alg:second-order}.
  \begin{algorithm}[H]
  \renewcommand{\thealgorithm}{4}
  \caption{Numerical algorithms for second-order case}
  \label{alg:second-order}
  \begin{algorithmic}[1]
    \Require The Brownian motion $ \Delta W(t_i) $, initial parameters $ (\theta^0,\tilde{p}_0^{0,\pi},\tilde{P}_0^{0,\pi}) $, learning rate $ \eta $;
    \Ensure Couple precess $ (\tilde{x}^{l,\pi}_{t_i},\tilde{u}^{l,\pi}_{t_i}) $ and $ \tilde{p}^{l,\pi}_{T} $.
    \For { $ l = 0 $ to $ maxstep $}
    \State $ \tilde{x}_{0}^{l,\pi} = x_{0} $, $ \tilde{p}_{0}^{l,\pi} = \tilde{p}_{0}^{l,\pi}, \tilde{P}_0^{l,\pi} = \tilde{P}_0^{l,\pi}; $
    \For { $ i = 0 $ to $ N-1 $}
    \State $\tilde{q}^{l,\pi}_{t_i} = \phi^1(t_i,\tilde{x}^{l,\pi}_{t_i};\theta^{l,1});$
    \State $\tilde{Q}^{l,\pi}_{t_i} = \phi^2(t_i,\tilde{x}^{l,\pi}_{t_i};\theta^{l,2});$
    \State $\tilde{u}^{l,\pi}_{t_i} = \arg\max_{u\in U}\mathcal{H}(t_i,\tilde{x}^{l,\pi}_{t_i},u);$
    \State $\tilde{x}^{l,\pi}_{t_{i+1}} = \tilde{x}^{l,\pi}_{t_i} + b(t_i,\tilde{x}^{l,\pi}_{t_i},\tilde{u}^{l,\pi}_{t_i})\Delta t_i + \sigma(t_i,\tilde{x}^{l,\pi}_{t_i},\tilde{u}^{l,\pi}_{t_i})\Delta W_{t_i};$
    \State $\tilde{p}^{l,\pi}_{t_{i+1}} = \tilde{p}^{l,\pi}_{t_i}-H_x(t_i,\tilde{x}^{l,\pi}_{t_i},\tilde{u}^{l,\pi}_{t_i},\tilde{p}^{l,\pi}_{t_i},\tilde{q}^{l,\pi}_{t_i})\Delta t_i + \tilde{q}^{l,\pi}_{t_i}\Delta W_{t_i};$
    \State $\tilde{P}^{l,\pi}_{t_{i+i}} = \tilde{P}^{l,\pi}_{t_i}- F(t_i,\tilde{x}^{l,\pi}_{t_i},\tilde{u}^{\pi}_{t_i},\tilde{p}^{l,\pi}_{t_i},\tilde{q}^{l,\pi}_{t_i},\tilde{P}^{l,\pi}_{t_i},\tilde{Q}^{l,\pi}_{t_i})\Delta t_i + \tilde{Q}^{l,\pi}_{t_i}\Delta W_{t_i};$
    \EndFor
    \State $ J(\tilde{u}^{l,\pi}(\cdot))=\dfrac{1}{M} \sum_{j=1}^{M}\Big[\dfrac{T}{N} \sum_{i=0}^{N-1}f(t_i,\tilde{x}_{t_i}^{l,\pi},\tilde{u}_{t_i}^{l,\pi}) + h(\tilde{x}_T^{l,\pi}) \Big]; $
    \State $ \mbox{loss} = \dfrac{1}{M} \sum_{j=1}^{M}\Big[ |-h_x(\tilde{x}^{l,\pi}_T)-\tilde{p}^{l,\pi}_T|^2 + |-h_{xx}(\tilde{x}^{l,\pi}_T)-\tilde{P}^{l,\pi}_T|^2 \Big]; $
    \State $ (\theta^{l+1},\tilde{p}_0^{l+1,\pi},\tilde{P}_0^{l+1,\pi})=(\theta^{l},\tilde{p}_0^{l,\pi},\tilde{P}_0^{l,\pi})-\eta\nabla\mbox{loss}. $
    \EndFor
  \end{algorithmic}
\end{algorithm}

Similar algorithm can be given when the state equation of a stochastic optimal control system is described by a fully coupled FBSDE.

\section{Existence and uniqueness results of FBSDEs}\label{appendix:Existence FBSDEs}
  For a special case of FBSDEs,
  \begin{equation}\label{eq:FBSDE_special}
    \begin{cases}
      X_{t} = X_{0} + \displaystyle \int_{0}^{t}b(s,X_{s},BY_{s},CZ_{s}), \mathrm{d} s + \int_{0}^{t}\sigma(s,X_{s},BY_{s},CZ_{s}) \mathrm{d} W_{s}, \vspace{1ex} \\
      Y_{t} = g(X_{T}) + \displaystyle \int_{t}^{T}f(s,X_{s},Y_{s},Z_{s})\, \mathrm{d} s - \int_{t}^{T} Z_{s}\mathrm{d} W_{s},
    \end{cases}
  \end{equation}
  where $ B $, $ C $ are $ k\times n $ matrixes, $ (x,y,z)\in\mathbb{R}^{n\times n\times n} $, and $ b,f,\sigma $ have appropriate dimensions.

  Denote that
  \begin{equation*}
    u = \left( \begin{array}{c} x \\ y \\z \end{array} \right), \ A(t, u)=\left(\begin{array}{c}
      -f \\b\\ \sigma
    \end{array}\right)(t,u),
  \end{equation*}
  and assume that
  \begin{assu}\label{assu:4}
    \begin{enumerate}
      \item $ A(t,u) $ is uniformly Lipschitz with respect to $ u $;
      \item $ A(\cdot,u) $ is in $ M^{2}(0,T) $ for $ \forall\ u $;
      \item $ g(x) $ is uniformly Lipschitz with respect to $ x \in \mathbb{R}^{n} $;
      \item $ g(x) $ is in $ L^{2}(\Omega,\mathcal{F}_{T},\mathbb{P}) $ for $ \forall\ x $;
      \item $ \forall x $, $ |l(t,x,By,Cz)-l(t,x,B\bar{y},C\bar{z})|\leq K(|B\hat{y}+C\hat{z}|) $, $ K>0 $, $ l=b,\sigma $
    \end{enumerate}
  \end{assu}
  and
  \begin{assu}\label{assu:5}
    \begin{align*}
      \left\langle A(t,u)-A(t,\bar{u}), u-\bar{u} \right\rangle & \leq - \nu_{1}|\hat{x}|^{2} - \nu_{2}|\hat{y}+\hat{z}|^{2},\\
      \left\langle g(x)-g(\bar{x}), (x-\bar{x}) \right\rangle & \geq 0,
    \end{align*}
    $$ \forall u =(x,y,z), \bar{u}=(\bar{x},\bar{y},\bar{z}), \hat{x} = x-\bar{x}, \hat{y} = y-\bar{y}, \hat{z} = z-\bar{z}, $$
    where $ \nu_{1} $ and $ \nu_{2} $ are given objective constants.
  \end{assu}

  The following result can be found in Theorem 2.6 of \cite{Peng1999Fully}.
  \begin{thm}\label{thm:special}
    Let Assumptions~\ref{assu:4} and \ref{assu:5} hold, then there exist a unique adapted solution $ (X,Y,Z) $ of FBSDE \eqref{eq:FBSDE_special}.
  \end{thm}

\bibliographystyle{ieeetr}
\bibliography{ref}

\begin{thebibliography}{10}

\bibitem{Bismut1972}
J.-M. Bismut, ``Analyse convexe et probabilitiés,'' {\em Thesis}, 1973.

\bibitem{Bismut1978}
J.-M. Bismut, ``An introductory approach to duality in optimal stochastic
  control,'' {\em SIAM Review}, vol.~20, no.~1, pp.~62--78, 1978.

\bibitem{Bensoussan1983Stochastic}
A.~Bensoussan, ``Stochastic maximum principle for distributed parameter
  system,'' {\em Journal of the Franklin Institute}, vol.~315, no.~5-6,
  pp.~387--406, 1983.

\bibitem{pontrygin1987}
L.~S. Pontrygin, ``Mathematical theory of optimal processes,'' {\em CRC Press},
  1987.

\bibitem{Bellman1958Dynamic}
Bellman and Richard, ``Dynamic programming and stochastic control processes,''
  {\em Information and Control}, vol.~1, no.~3, pp.~228--239, 1958.

\bibitem{Harold_numerical}
H.~Kushner and P.~G. Dupuis, {\em Numerical Methods for Stochastic Control
  Problems in Continuous Time}.
\newblock Springer, 2001.

\bibitem{Quadrat1994Numerical}
H.~J. Kushner, ``Numerical methods for stochastic control problems in
  continuous time,'' {\em SIAM, J. Control Optim.}, vol.~28, pp.~888--1048,
  1990.

\bibitem{dong2007the}
H.~Dong and N.~V. Krylov, ``The rate of convergence of finite-difference
  approximations for parabolic bellman equations with lipschitz coefficients in
  cylindrical domains,'' {\em Applied Mathematics and Optimization}, vol.~56,
  no.~1, pp.~37--66, 2007.

\bibitem{jakobsen2003on}
E.~R. Jakobsen, ``On the rate of convergence of approximation schemes for
  bellman equations associated with optimal stopping time problems,'' {\em
  Mathematical Models and Methods in Applied Sciences}, vol.~13, no.~05,
  pp.~613--644, 2003.

\bibitem{krylov2005the}
N.~V. Krylov, ``The rate of convergence of finite-difference approximations for
  bellman equations with lipschitz coefficients,'' {\em Applied Mathematics and
  Optimization}, vol.~52, no.~3, pp.~365--399, 2005.

\bibitem{bertsekas1995neuro-dynamic}
D.~P. Bertsekas and J.~N. Tsitsiklis, ``Neuro-dynamic programming: an
  overview,'' {\em Proceedings of 1995 34th IEEE Conference on Decision and
  Control}, vol.~1, pp.~560--564, 1995.

\bibitem{PardalosApproximate}
Pardalos and M.~Panos, ``Approximate dynamic programming: solving the curses of
  dimensionality,'' {\em Optimization Methods and Software}, vol.~24, no.~1,
  pp.~155--155, 2009.

\bibitem{BengioDL}
I.~{Goodfellow}, Y.~{Bengio}, and A.~{Courville}, {\em Deep Learning}.
\newblock 2016.

\bibitem{han2016deep}
J.~Han and W.~E, ``Deep learning approximation for stochastic control
  problems,'' {\em Deep Reinforcement Learning Workshop}, 2016.

\bibitem{WeinanDLforBSDE}
W.~E, J.~Han, and A.~Jentzen, ``Deep learning-based numerical methods for
  high-dimensional parabolic partial differential equations and backward
  stochastic differential equations,'' {\em Communications in Mathematics and
  Statistics}, vol.~5, no.~4, pp.~349--380, 2017.

\bibitem{HanPNAS}
J.~Han, A.~Jentzen, and W.~E, ``Solving high-dimensional partial differential
  equations using deep learning,'' {\em Proceedings of the National Academy of
  Sciences}, vol.~115, no.~34, pp.~8505--8510, 2018.

\bibitem{deeplearning_FBSDE}
J.~Han and J.~Long, ``Convergence of the deep bsde method for coupled fbsdes,''
  {\em arXiv:1811.01165}, 2018.

\bibitem{Peng_FBSDE_numerical}
S.~Ji, S.~Peng, Y.~Peng, and X.~Zhang, ``Three algorithms for solving
  high-dimensional fully coupled fbsdes through deep learning,'' {\em IEEE
  Intelligent Systems}, vol.~35, no.~3, pp.~71--84, 2020.

\bibitem{hure2020deep}
C.~Hur{\'e}, H.~Pham, and X.~Warin, ``Deep backward schemes for
  high-dimensional nonlinear pdes,'' {\em Mathematics of Computation}, vol.~89,
  no.~324, pp.~1547--1579, 2020.

\bibitem{raissi2018forward-backward}
M.~Raissi, ``Forward-backward stochastic neural networks: Deep learning of
  high-dimensional partial differential equations,'' {\em arXiv: 1804.07010},
  2018.

\bibitem{pham2018deep_1}
C.~Huré, H.~Pham, A.~Bachouch, and N.~Langrené, ``Deep neural networks
  algorithms for stochastic control problems on finite horizon: Convergence
  analysis,'' {\em SIAM Journal on Numerical Analysis}, vol.~59, no.~1,
  pp.~525--557, 2021.

\bibitem{pham2018deep_2}
A.~Bachouch, C.~Huré, N.~Langrené, and H.~Pham, ``Deep {Neural} {Networks}
  {Algorithms} for {Stochastic} {Control} {Problems} on {Finite} {Horizon}:
  {Numerical} {Applications},'' {\em Methodology and Computing in Applied
  Probability}, 2021.

\bibitem{pereira2019neural}
M.~A. Pereira, Z.~Wang, I.~Exarchos, and E.~A. Theodorou, ``Neural network
  architectures for stochastic control using the nonlinear feynman-kac lemma,''
  {\em arXiv:1902.03986v2}, 2019.

\bibitem{germain2021neural}
M.~Germain, H.~Pham, and X.~Warin, ``Neural networks-based algorithms for
  stochastic control and pdes in finance,'' {\em arXiv preprint
  arXiv:2101.08068}, 2021.

\bibitem{Peng90}
S.~Peng, ``A general stochastic maximum principle for optimal control
  problems,'' {\em Siam Journal on Control and Optimization}, vol.~28, no.~4,
  pp.~966--979, 1990.

\bibitem{Yong_stochastic_control}
J.~Yong and X.~Zhou, {\em Stochastic Controls-Hamiltonian System and HJB
  Equations}.
\newblock Springer, 1999.

\bibitem{Bensoussan1982}
A.~Bensoussan, {\em Lecture on stochastic control in Nonlinear Filtering and
  Stochastic Control}.
\newblock Springer, 1982.

\bibitem{LMFGS1989}
D.~C. Liu and J.~Nocedal, ``On the limited memory bfgs method for large scale
  optimization,'' {\em Mathematical Programming}, vol.~45, no.~1-3,
  pp.~503--528, 1989.

\bibitem{Peng1999Fully}
S.~Peng and Z.~Wu, ``Fully coupled forward-backward stochastic differential
  equations and applications to optimal control,'' {\em Siam Journal on Control
  and Optimization}, vol.~37, no.~3, pp.~825--843, 1999.

\bibitem{ma1999forward}
J.~Ma and J.~Yong, {\em Forward-backward stochastic differential equations and
  their applications}.
\newblock No.~1702, Springer Science \& Business Media, 1999.

\bibitem{Peng2000Problem}
S.~Peng, ``Problem of eigenvalues of stochastic hamiltonian systems with
  boundary conditions,'' {\em Stochastic Processes \& Their Applications},
  vol.~88, no.~2, pp.~259--290, 2000.

\bibitem{Peng1993Backward}
S.~Peng, ``Backward stochastic differential equations and applications to
  optimal control,'' {\em Applied Mathematics and Optimization}, vol.~27,
  no.~2, pp.~125--144, 1993.

\bibitem{hu2018a}
M.~Hu, S.~Ji, and X.~Xue, ``A global stochastic maximum principle for fully
  coupled forward-backward stochastic systems,'' {\em Siam Journal on Control
  and Optimization}, vol.~56, no.~6, pp.~4309--4335, 2018.

\bibitem{peng2019nonlinear}
S.~Peng, {\em Nonlinear Expectations and Stochastic Calculus under
  Uncertainty-with Robust CLT and G-Brownian Motion}.
\newblock Springer, 2019.

\end{thebibliography}

\end{document}